\documentclass{amsart}[11pt]

\usepackage[usenames]{color}
\usepackage{amssymb}
\usepackage{graphicx, epsfig}
\usepackage{latexsym, amsfonts, amscd, amsmath}
\usepackage{mathrsfs}
\makeindex 
\input xy
\xyoption{all}

\definecolor{Indigo}{rgb}{0.2,0.1,0.7}
\definecolor{Violet}{rgb}{0.5,0.1,0.7}

\newtheorem{thm}{Theorem}[section]
\newtheorem{prop}[thm]{Proposition}
\newtheorem{lemma}[thm]{Lemma}
\newtheorem{cor}[thm]{Corollary}

\theoremstyle{definition}
\newtheorem{defn}[thm]{Definition}

\theoremstyle{remark}
\newtheorem{rem}[thm]{Remark}

\newenvironment{mat}{\left(\begin{smallmatrix}
\end{smallmatrix}\right)}{enddef}

\numberwithin{equation}{section} \numberwithin{figure}{section}
\numberwithin{table}{section}

\newcommand{\End}{{\operatorname{End}}}

\newcommand{\Ker}{{\operatorname{Ker}}}

\newcommand{\Spec}{{\operatorname{Spec }}}
\newcommand{\Spf}{{\operatorname{Spf }}}
\newcommand{\Sp}{{\operatorname{Sp }}}

\newcommand{\val}{{\operatorname{val}}}

\newcommand{\GL}{{\operatorname{GL}}}

\newcommand{\Lie}{{\operatorname{Lie}}}

\newcommand{\gerh}{{\frak{h}}}

\newcommand{\germ}{{\frak{m}}}

\newcommand{\gers}{{\frak{s}}}
\newcommand{\gert}{{\frak{t}}}

\newcommand{\gerF}{{\frak{F}}}

\newcommand{\gerN}{{\frak{N}}}

\newcommand{\gerX}{{\frak{X}}}
\newcommand{\gerY}{{\frak{Y}}}
\newcommand{\gerZ}{{\frak{Z}}}

\newcommand{\calA}{{\mathcal{A}}}

\newcommand{\calD}{{\mathcal{D}}}

\newcommand{\calF}{{\mathcal{F}}}
\newcommand{\calG}{{\mathcal{G}}}

\newcommand{\calO}{{\mathcal{O}}}
\newcommand{\calP}{{\mathcal{P}}}

\newcommand{\calU}{{\mathcal{U}}}
\newcommand{\calV}{{\mathcal{V}}}
\newcommand{\calW}{{\mathcal{W}}}
\newcommand{\calX}{{\mathcal{X}}}
\newcommand{\calY}{{\mathcal{Y}}}
\newcommand{\calZ}{{\mathcal{Z}}}

\def\AA{\mathbb{A}}

\def\CC{\mathbb{C}}

\def\FF{\mathbb{F}}

\def\QQ{\mathbb{Q}}
\def\RR{\mathbb{R}}

\def\ZZ{\mathbb{Z}}

\newcommand{\scrP}{{\mathscr{P}}}

\newcommand{\id}{{\noindent}}

\newcommand{\Yrig}{{\gerY_\text{\rm rig}}}

\newcommand{\pirig}{{\pi_\text{\rm rig}}}
\newcommand{\rig}{{\operatorname{rig}}}

\newcommand{\pr}{{\rm pr}}

\newcommand{\ra}{{\rightarrow}}
\newcommand{\inc}{\hookrightarrow}
\newcommand{\lra}{\longrightarrow}              

\newcommand{\llra}[1]{\stackrel{#1}{\lra}}      
\newcommand{\eqra}{\llra{\sim}}                 
\newcommand{\U}{{\rm U}}
\newcommand{\Xm}{X(m)}
\newcommand{\Xmo}{X^0(m)}
\newcommand{\gerXm}{\gerX(m)}
\newcommand{\gerXmo}{\gerX^0(m)}

\newcommand{\ts}{\otimes}
\newcommand{\tsc}{\hat{\ts}}                    
\newcommand{\omg}{\underline{\omega}}
\begin{document}


\title{Overconvergence and classicality: the case of curves}
\author{Payman L Kassaei}
\email{kassaei@alum.mit.edu}
\address{Department of Mathematics, King's College London, Strand, London, WC2R 2LS, UK. }

\maketitle

\begin{abstract}
Given our set up of  a system of abstract curves and maps between them satisfying certain assumptions, we  prove a classicality criterion for overconvergent sections of line bundles over these curves.
As a result we prove such criteria for overconvergent modular forms over various Shimura curves. In particular we provide a classicality criterion for overconvergent modular forms studied in \cite{Kassaei2} and their higher-level generalizations.

\end{abstract}
\section{Introduction}

Over a decade has passed since R. Coleman's breakthrough in applying the theory of overconvergent modular forms to the study of $p$-adic variation of modular forms of finite slope.
Coleman's results have been extended and generalized in various directions and have found applications, most notably via constructing $p$-adic families of automorphic forms. One ingredient, however, is still missing in some important generalizations
of Coleman's theory, mainly in those which use the  geometry of Shimura varieties in the construction and study of overconvergent automorphic forms. That ingredient is Coleman's
classicality result \cite{Coleman1,Coleman2} which states that ``overconvergent modular forms of small slope are classical'' and which is crucial  in most applications of the
theory; for example,  to construct a $p$-adic family of  classical automorphic forms containing a given one, one uses methods of $p$-adic analysis to first construct  a family of overconvergent automorphic forms (objects more inviting of $p$-adic interpolation) and then one invokes a classicality result to deduce that all but finitely many members of the family are indeed classical automorphic forms.   Somehow, it has not been easy  to extend Coleman's clever dimension-counting proof of the classicality result  to other cases. The ordinary case (i.e., when slope is zero) was dealt with by
Hida and has been extensively generalized by Hida and others.

In \cite{Kassaei3} we presented an alternative proof of Coleman's result which was based on the formal and rigid geometry of the modular curves.
The strategy is to p-adically analytically continue an overconvergent modular form to a global p-adic analytic section of a line bundle which will then, according to a rigid geometric GAGA, be   a classical (algebraic) modular form. The analytic continuation consists of two steps: first we use Buzzard's analytic continuation results \cite{Buzzard} to extend an overconvergent modular form to the entire supersingular locus, and then, we construct a second modular form on the complementary region and show that the two can be ``glued" together despite the fact that their regions of definitions are disjoint. A good part of the work goes into the gluing process and uses the full force of the classical theory of canonical subgroups of elliptic curves.

In this paper we show that this method can be applied in the context of various Shimura curves. In particular, our results provide a classicality result for the overconvergent modular forms over unitary Shimura curves  which were studied in \cite{Kassaei2}. In fact we generalize the basic constructions of that paper to the case where the level has arbitrary powers of $\calP$ in it, and prove the classicality result in that generality. We also prove a classicality result for the quaternionic overconvergent modular forms studied in \cite{Kassaei1}.

Our original presentation of the proof in \cite{Kassaei3} relies (seemingly) essentially on the moduli problems that the modular curves in question are a solution to. As we were trying to carry out this method over Shimura curves it became clear that the existence of a moduli problem is a bit of a red herring. It is instead  some specific formal and rigid geometric features of the Shimura curves and of certain maps between them that are at work. In this article we take this viewpoint and show that given a system of abstract curves and maps between them satisfying certain assumptions, one can develop a theory of overconvergent sections of line bundles on these curves and prove a classicality criterion. Our assumptions are general enough to cover all  cases of Shimura curves that we are interested in.
One reason that the argument can work in this generality is an equally general treatment of all desired aspects of a theory of canonical subgroups for curves in \cite{GorKas}.

We are working on proving similar results in higher dimensions in a similar spirit. Some of the constructions in this paper can be carried out in dimensions bigger than one, but in general there are a few obstacles in proving a classicality result. To begin with one needs a comprehensive enough theory of canonical subgroups (which is currently the subject of research of a number of people).  Recently Shu Sasaki has used the above analytic continuation method to prove a classicality criterion for overconvergent  Hilbert modular forms when $p$ is split in the totally real field in question. In that case canonical subgroups can be constructed using the classical method of Lubin-Katz, as when $p$ is split the formal group of  an HBAV  factorizes as a product of one-dimensional formal groups. When $p$ is not split, canonical subgroups are not well understood yet. Furthermore, in this case, Buzzard's analytic continuation method does not automatically extend the overconvergent  Hilbert  modular form to the entire non-ordinary locus of the Hilbert modular variety.

The article has four parts. In \S 2 we introduce our set-up and define spaces of overconvergent sections of line bundles and define the action of a completely continuous $\U$ operator. In \S 3 we
use Buzzard's method of successive hitting by the $\U$ operator to provide partial analytic continuation of overconvergent sections.
In \S 4 we carry out our method of analytic continuation and prove the classicality criterion. In \S 5 we show how these general results can be applied over Shimura curves.

\

\noindent {\bf{Acknowledgments:}}  I am grateful to the Max-Planck-Institut f\"{u}r Mathematik in Bonn for their hospitality and the excellent work conditions provided to me in the Spring of 2005 when a major part of this work was carried out.
 \section{Overconvergence}
\subsection{Set-up: the ``tame'' case}

 \id Let~$p$ be a prime and $L_0$ a finite extension of $\QQ_p$. Let $\calO_0$ denote the ring of integers with maximal ideal $\germ=(\varpi)$, and residue field
 $\kappa\cong\FF_q$. Let $\val$ be a valuation normalized so that
 $\val(\varpi)=1$.  Define $|.|=|.|_{L_0} =(1/q)^{\val(.)}$ on $L_0$. Let $\hat{\bar{L}}_0$
denote the completion of a fixed algebraic closure of $L_0$. The
valuation on $L_0$ can be extended to a valuation
$\val:\hat{\bar{L}}_0 \ra \QQ$, and hence the absolute value $|.|$
too can be extended to $\hat{\bar{L}}_0$.  If $L$ is a completely
valued subfield of $\hat{\bar{L}}_0$, we define $\val$ and $|.|$ on
$L$ by restriction from $\hat{\bar{L}}_0$.

Let~$R$ be an~$\calO_0$-algebra. By a ``curve"~$Z$ over~$R$ we mean
a flat finite-type separated morphism~$f\colon Z \rightarrow
\Spec(R)$   such that the geometric fibres of~$f$ are connected and
of dimension one, and that $Z$ is a reduced  scheme. If $S$ is a
scheme over $R$, by $Z \ts S$ we mean the base change of $Z$ via $S
\ra \Spec(R)$. If $S=\Spec(R^\prime)$ we denote $Z \ts S$ also by $Z
\ts {R^\prime}$. This convention applies in the same way to all
other relative objects in this paper. However, we often denote a
morphism and its base change by the same notation.

Guided by examples of Shimura curves, we introduce the following
data. Let~$X,Y$ be curves over~$\calO_0$ with a morphism $\pi:Y \ra X$ such that

\begin{itemize}
\item[\bf{A1}\ ] $X$ is smooth over $\calO_0$;

\item[\bf{A2}\ ]$Y$ is a regular scheme such that

\begin{itemize}
 \item[\bf{A2.1}\ ] there exists a section~$s\colon X \ts \kappa \ra
  Y \ts \kappa$ to~$\pi \ts \kappa:Y\ts \kappa \ra X \ts \kappa$,

  \item[\bf{A2.2}\ ] the special fibre~$Y \ts \kappa$ is reduced, has two
  components, and each intersection point  of the components is defined over~$\kappa$ and
its completed local ring is isomorphic to~$\kappa[\![s, t]\!]/(st)$.

  \item[\bf{A2.3}\ ] the set theoretic
preimage~$(\pi \ts \kappa)^{-1}(\pi \ts \kappa)(Q)$ is equal
to~$Q$ for any singular point~$Q\in Y \ts \kappa$,

  \item[\bf{A2.4}\ ] we have fixed an automorphism $w:Y \ra Y$ defined over $\calO_0$ whose
reduction mod $\varpi$ switches the  components of $Y \ts
\kappa$,

\item[\bf{A2.5}\ ] we have fixed an automorphism $\delta: Y \ra Y$ defined over $\calO_0$ whose reduction mod $\varpi$  sends each component of $Y \ts \kappa$ to itself;

\end{itemize}

\item[\bf{A3}\ ] the morphism  $\pi:Y \ra X$ is finite flat of
degree $1+e$ where $e>1$ is an integer.
\end{itemize}
We define~$(Y \ts\kappa)^{(\infty)}=s(X \ts \kappa)\setminus (Y \ts
\kappa)^{\rm sing}$, and~$(Y \ts \kappa)^{(0)}=(Y \ts \kappa)
\setminus s(X \ts \kappa)$.

\begin{rem}{\label{finite extension}}
Assumptiuon {\bf A2.2} can be relaxed. It is enough to assume that
every point of intersection is defined over a finite extension
$\kappa^\prime$ of $\kappa$ and that $Y \ts \kappa^\prime$ is
reduced, has two components, and the completed local ring of each
intersection point of those components is isomorphic to
$\kappa^\prime[\![s,t]\!]/(st)$.
\end{rem}
\begin{rem}{\label{analogy-levelN}}
We will apply the results of this paper to examples where $X,Y$ are
various types of Shimura curves. As an  example of the above
consider $X$ to be a modular curve with level prime to $p$, say
$X(\Gamma_1(N))$ with $(p,N)=1$, whose noncuspidal points classify
$(E,i)$ with $E$ an elliptic curve, and $i$ a $\Gamma_1(N)$-level
structure.  Also let $Y$ be obtained from $X$ by adding a
$\Gamma_0(p)$-level structure, say $X(\Gamma_1(N) \cap
\Gamma_0(p))$, whose noncuspidal points classify $(E,i,C)$, where
$(E,i)$ is as above and $C$ is a subgroup of order $p$ in $E$. The
automorphism $w$ can be then given by dividing out an elliptic curve
and its  $\Gamma_1(N)$-level structure by its $\Gamma_0(p)$-level
structure, and adding the  $\Gamma_0(p)$-level structure induced by
the $p$-torsion points after passage to the quotient. The
automorphism $\delta$ is the diamond operator $<p^{-1}>$. In this
case $e=p$ and the section $s$ is defined by adding the
$\Gamma_0(p)$-level structure given by the subgroup scheme
$\Ker({\rm Frob}_p)$ in characteristic $p$.

\end{rem}

Let~$\gerX, \gerY$ be the formal schemes obtained, respectively,
by completing~$X, Y$ along their special fibres. These are
quasi-compact admissible formal schemes over $\calO_0$ as defined
in \cite{BoschLutkebohmertI}. For simplicity, we use the same
notation for the induced morphisms after formal completion. There
is a functor ``$\rig$" which associates to every quasi-compact
admissible formal scheme $\gerZ$ over $\calO_0$, its ``generic
fibre", $\gerZ_\rig$, which is a quasi-compact and quasi-separated
rigid analytic space over $L_0$. See \cite{BoschLutkebohmertI} for
an account of this construction which is due to Raynaud, or \S 2.1
of \cite{GorKas} for a brief survey. We denote the image of a map
$\alpha$ under this functor by $\alpha_\rig$.

Since we will frequently use results from \cite{GorKas}, the
notation has been chosen in accordance with that article, except
that $\calO_0,L_0$ are denoted by $\calO,K$ there. In \S 2.3 of {\it
loc. cit.} a ``measure of singularity" is defined which is modeled
over the notion of the measure of supersingularity of elliptic
curves in the context of modular curves. For a point $P$ of
$\gerX_\rig$ we have~$\nu_\gerX(P)\in \QQ^{\geq 0}$, which is
well-defined only when $\nu_\gerX(P)<1$ (and the statement
``$\nu_\gerX(P)\geq1$''  is also well-defined). For a point $Q \in
\gerY_\rig$ we have ~$\nu_\gerY(Q)\in \QQ^{\geq 0}$ (always well
defined, and at most $1$). Over a residue annulus in $\gerY_\rig$,
$\nu_\gerY$ is the valuation of a carefully chosen parameter. For a
point $Q$ outside the union of the residue annuli of singular points
of $Y \ts \kappa$, one has $\nu_\gerY(Q)=0$ or $1$, depending on
whether $Q$ specializes to a point in $(Y \ts \kappa)^{(\infty)}$ or
$(Y \ts \kappa)^{(0)}$, respectively. We refer to \S 2.3 of {\it
loc. cit.} for precise definitions. For every interval~$I \subset
\RR$ with endpoints in $\QQ$ there is an admissible open~$\gerY_\rig
I$ in $\gerY_\rig$ whose points are~
\[
\{Q\in \Yrig: \nu_\gerY(Q) \in I\}.
\]
A similar notation will be used for $\gerX_\rig$. In this case,
however, the interval $I$ is assumed to be inside $[0,1)$. If $L
\subset \hat{\bar{L}}_0$ is a completely valued extension of $L_0$
then $\nu_\gerX$ and $\nu_\gerY$ can be defined over $\gerX_\rig
\tsc L$ and $\gerY_\rig \tsc L$ by pullback.

In \S 3 of \cite{GorKas} it is
proven that the morphism $\pi_\rig:\gerY_\rig \ra \gerX_\rig$
admits a section
\[
\gers_\rig:\gerX_\rig[0,e/(e+1)) \ra \gerY_\rig[0,e/(e+1))
\]
which we call the {\it canonical} section. To simplify the
notation, for a point $Q \in \gerY_\rig$, we sometimes denote
$w_\rig(Q)$ by $Q^w$. We also denote the base extension of $\gers_\rig$ to any extension $L$ of $L_0$ by the same notation.

We summarize some results of \cite{GorKas} in the following
proposition. We refer to Definition 3.11 and Lemmas 3.6 and 4.2  of
\cite{GorKas} for details. Item (5) does not appear in \cite{GorKas} but can be proven in the same way as item (4). Note that all these
results are proven for $L=L_0$ in \cite{GorKas} but the results for general $L$ follow immediately.

\begin{prop}{\label{proposition:GK}}
Let $L\subset \hat{\bar{L}}_0$ be a completely valued extension of $L_0$. Let $Q$ be a point of $\gerY_\rig \tsc L$.
\begin{itemize}
\item[(1)] If $\nu_\gerY(Q)<e/(e+1)$, then
$\nu_\gerY(Q)=\nu_\gerX(\pirig Q)$. In this case, we say~$Q$ is
canonical. A point $Q$ is canonical if and only if it is in the
image of the canonical section $\gers_\rig$.

\item[(2)]  If~$\nu_\gerY(Q) > e/(e+1)$, then we have
$\nu_\gerY(Q)=1-e^{-1}\nu_\gerX(\pirig Q)$. In this case, we say
that~$Q$ is {\it anti-canonical}.

\item[(3)]  We say that~$Q$ is {\it too singular} if
$\nu_\gerY(Q)=e/(e+1)$. This is equivalent to~$\nu_\gerX(\pirig Q)
\geq e/(e+1)$.

\item [(4)] We have $\nu_\gerY(w_\rig(Q))=1-\nu_\gerY(Q)$.

\item [(5)] We have $\nu_\gerY(\delta_\rig(Q))=\nu_\gerY(Q)$.
\end{itemize}

\end{prop}

\begin{cor}{\label{corollary:positions}}
Let $L\subset \hat{\bar{L}}_0$ be a completely valued extension of $L_0$. Let $Q_1$ and $Q_2$ be points on $\gerY_\rig \tsc L$ so that $\pi_\rig(Q_1)=\pi_\rig(Q_2)=P \in \gerX_\rig \tsc L$. Then

\begin{itemize}
\item[(1)] The point $Q_1$ is too singular iff $Q_2$ is too singular iff $\nu_\gerX(P) \geq e/(e+1)$. In that case $\nu_\gerY(Q_2^w)=1/(e+1)$.

\item[(2)] If $Q_1$ and $Q_2$ are both canonical, or both anti-canonical, then $\nu_\gerY(Q_2^w)=1-\nu_\gerY(Q_1)$.

\item[(3)] If $Q_1$ is canonical and $Q_2$ is anti-canonical, then
$\nu_\gerY(Q_2^w)=e^{-1}\nu_\gerY(Q_1)$.

\item[(4)] If $Q_1$ is anti-canonical and $Q_2$ is canonical, then
$\nu_\gerY(Q_2^w)=1-e(1-\nu_\gerY(Q_1))$.

\end{itemize}
\end{cor}

\begin{proof}
All the statements follow easily from Proposition
\ref{proposition:GK}. Note that in (2) if both $Q_1$ and $Q_2$ are
canonical they must be equal.
\end{proof}

\begin{lemma}{\label{lemma:affinoid}}
If $I\subsetneq [0,1]$ (respectively, $I \subset [0,1)$) is a closed
interval with endpoints in $\QQ$, then $\gerY_\rig I$ (respectively,
$\gerX_\rig I$) is an  affinoid subdomain.
\end{lemma}
\vspace{-2mm}

\begin{proof}
Any finite union of affinoids on an irreducible curve is either
the whole space or an affinoid itself. It is therefore enough to show
$\gerX_\rig I$ and $\gerY_\rig I$ are quasi-compact subdomains.

Over $\gerX_\rig$  there is another well-known general construction
which will produce the domains $\gerX_\rig[a,b]$. See, for example,
\S 3.2 of \cite{KisinLai}. In the notation of that paper, let
$\calD$ be the Cartier divisor on $X \ts \kappa$ given by the sum of
all points $(\pi \ts \kappa) (\beta)$ where $\beta$ runs over
singular points of $Y \ts \kappa$. It is easy to show that
$\gerX_\rig[0,r]$  in our construction is the same as
$\gerX_\rig(p^{-r})$ defined there (which is evidently
quasi-compact), and in fact this can be done for any closed interval
. One point to remember is that based on our choice of the valuation
the ramification degree $e$ considered in \cite{KisinLai} is equal
to 1 here.

Over $\gerY_\rig$, however, the above construction doesn't work. We
refer the reader to \S 2.3 of \cite{GorKas} for details on the
definition of $\nu_\gerY$.   The domain $\gerY_\rig(0,1)$ is a
finite disjoint union of open annuli $\{x: 0 < v(x)<1\}$ where $x$
is the specific (type of)  parameter used in the definition of
$\nu_\gerY$. Therefore, if $[a,b] = I \subset (0,1)$ we can think of
$\gerY_\rig I$ as an admissible finite disjoint union of closed
annuli $\{x: a \leq |x| \leq b\}$.  This shows that $\gerY_\rig I$
is an affinoid subdomain in this case. Finally to address cases
where exactly one of $0$ or $1$ belongs to $I$ we note that it is
enough to consider intervals of the form $[0,r]$ (respectively
$[r,1]$) where $r<e/(e+1)$ (respectively $r>e/(e+1)$). The reason is
that, for example, if $r\geq e/(e+1)$ then
$\gerY_\rig[0,r]=\gerY_\rig[0,1/(e+1)] \cup \gerY_\rig[1/(e+1),r]$
which is quasi-compact since $\gerY_\rig$ is quasi-separated and
therefore any finite covering by quasi-compact opens is an
admissible covering.  If $r<e/(e+1)$ then $\gerY_\rig[0,r]$ is
isomorphic to $\gerX_\rig[0,r]$ by the existence of $\gers_\rig$. If
$r>e/(e+1)$, then $\gerY_\rig[r,1]$ is a connected component in
$\pi_\rig^{-1}(\gerX_\rig[0,e(1-r)])$ by Proposition
\ref{proposition:GK}. Since $\pi_\rig$ is a finite  morphism this
implies that $\gerY_\rig[r,1]$ is quasi-compact.
\end{proof}

We now introduce a curve over $\gerY_\rig$ which allows us to define
a correspondence on $\gerY_\rig$. See Remark \ref{analogy-Y^0} for
the analogue in the case of modular curves. Let $\gerY^0_\rig$ be a
rigid analytic curve over $\Spec(L_0)$ fitting into the following
 diagram. (The notation $\gerY^0_\rig$ is chosen in accordance with
the rest of the notation and is {\underline{not}} meant to suggest
that $\gerY^0_\rig$ is obtained via the process of formal completion
followed by applying the functor rig from some specific curve
``$Y^0$''. The same warning goes for the maps).

\begin{eqnarray}\label{diagram:tame}
\xymatrix{ \gerY^0_\rig\ar@{}@<.5ex>[r]^{\pi_{1,\rig}^\prime} \ar@<1ex>[r]
\ar@{}@<-.5ex>[r]_{\pi_{1,\rig}} \ar@<-1ex>[r] &\gerY_\rig \ar[d]^{\pi_\rig}\\
               & \gerX_\rig}
\end{eqnarray}
such that
\begin{itemize}

\item[\bf{A4}\ ]$\pi_{1,\rig}$, $\pi_{1,\rig}^\prime$ are finite flat rigid analytic morphisms
defined over $L_0$. Define $\pi_{2,\rig}:\gerY^0_\rig \ra \gerY_\rig$ by
$\pi_{2,\rig}:=w_\rig\pi_{1,\rig}^\prime$.

\item[\bf{A5}\ ]For any $Q \in \gerY_\rig^0$ we have $\pi_{1,\rig}(Q) \neq \pi^\prime_{1,\rig}(Q)$. In particular, for such a point $Q$ if one of  $\pi_{1,\rig}(Q)$  or $\pi_{1,\rig}^\prime(Q)$ is canonical, the other will be anti-canonical.

\end{itemize}
By Assumption {\bf A5}, we have
\begin{eqnarray}\label{equation:can-antican}
(\pi^\prime_{1,\rig})^{-1}(\gerY_\rig[0,e/(e+1))) \subset
(\pi_{1,\rig})^{-1}(\gerY_\rig(e/(e+1),1]). \
\end{eqnarray}
We further assume the following.
\begin{itemize}

\item[\bf{A6}\ ] By Diagram \ref{diagram:tame}, Equation
\ref{equation:can-antican}, and parts (1) and (2) of Proposition
\ref{proposition:GK} we have a commutative diagram

\begin{eqnarray}\label{diagram:product-tame}
\xymatrix{ (\pi^\prime_{1,\rig})^{-1}(\gerY_\rig[0,e/(e+1)))
\ar[rr]^{\pi_{1,\rig}} \ar[d]_{\pi^\prime_{1,\rig}} &&
\gerY_\rig(e/(e+1),1] \ar[d]^{\pi_\rig}  \\
\gerY_\rig[0,e/(e+1)) \ar[rr]_{\pi_\rig} && \gerX_\rig[0,e/(e+1))}
\end{eqnarray}
which we assume to be a {\underline{product}} diagram. In
particular, by base extension of $\gers_\rig$, there is a section
\[
\gers_\rig^0:\gerY_\rig(e/(e+1),1]  \ra \gerY_\rig^0
\]
to $\pi_{1,\rig}$ whose image is
$(\pi^\prime_{1,\rig})^{-1}(\gerY_\rig[0,e/(e+1)))$. This is simply equivalent to requiring that the top arrow is an isomorphism and $\gers_\rig^0$ is its inverse.

\end{itemize}

\begin{rem}\label{analogy-Y^0}
In applications to Shimura curves, $\gerY^0_\rig$ will be the analogue of the modular curve $X(\Gamma_1(N) \cap \Gamma^0(p))$ whose noncuspidal points classify  $(E,i,C,D)$  with $E$ an elliptic curve, $i$ a  $\Gamma_1(N)$-level structure, and $C,D$   two finite-flat subgroups of order $p$ which intersect trivially.  The morphisms $\pi_{1,\rig}$ and $\pi_{1,\rig}^\prime$  are then the ones that forget $D$ and   $C$, respectively. The morphism $\pi_{2,\rig}$ is the one that quotients out by $D$. In the classical theory this data is used to define the Hecke correspondence
$\U_p$.
\end{rem}

Here we prove a Lemma that we will use many times in this paper.
For an interval $I$ define
\begin{eqnarray}
I^{\tau}&=&\{1-e(1-t)| t \in I\} \\
I^{w}&=&\{1-t| t \in I\}\\
 I^{\sigma}&=&\{e^{-1}t| t\in I\}.
\end{eqnarray}

\begin{lemma}{\label{lemma:pre-trace}}
\begin{itemize}
\item[(1)] If $I \subseteq [0,e/(e+1))$, i.e., inside the canonical locus we have
\[
\pi_{1,\rig}^{-1}(\gerY_\rig I)\subseteq
\pi_{2,\rig}^{-1}(\gerY_\rig I^{\sigma}).
\]

\item[(2)] Over the too singular locus we have
\[
\pi_{1,\rig}^{-1}(\gerY_\rig [e/(e+1),e/(e+1)])\subseteq
\pi_{2,\rig}^{-1}(\gerY_\rig [1/(e+1),1/(e+1)]).
\]

\item[(3)] If $I \subseteq (e/(e+1),1]$, i.e., inside the anti-canonical locus we have
\[
\pi_{1,\rig}^{-1}(\gerY_\rig I)\subseteq
\pi_{2,\rig}^{-1}(\gerY_\rig  I^{w}) \cup \pi_{2,\rig}^{-1}(\gerY_\rig I^{\tau})
\]
where the right hand side is an admissible disjoint union.

\end{itemize}
\end{lemma}

\begin{proof}
Let $Q \in  \pi_{1,\rig}^{-1}(\gerY_\rig I)$. Define $Q_{1}:= \pi_{1,\rig}(Q)$ and $Q_{2}:= \pi^{\prime}_{1,\rig}(Q)$. Then $\pi_{2,\rig}(Q)=Q_{2}^{w}$. Also from Diagram \ref{diagram:tame} we have $\pi_{\rig}(Q_{1})=\pi_{\rig}(Q_{2})$ and hence Corollary \ref{corollary:positions} can be applied.

In case (1), we know $Q_{1}$ is canonical  and thus assumption {\bf A5}
tells us that $Q_{2}$ is anti-canonical. Hence the result follows from part (3) of Corollary \ref{corollary:positions}. In case (2), we know $Q_{1}$ is too singular and hence the result follows from part (1) of Corollary \ref{corollary:positions}. Finally, in case (3), the point $Q_{1}$ is anti-canonical and parts (2) and (4) of Corollary \ref{corollary:positions} imply the result. For the last statement
 notice that under the assumptions on $I$ we have $I^{w} \cap I^{\tau}=\emptyset$, and hence the union is disjoint. If $I$ is a closed interval with rational endpoints then using Lemma \ref{lemma:affinoid} it is clear that the covering is admissible (all our spaces are quasi-separated). The general case is reduced to this case using the maximum modulus principle.

\end{proof}

\begin{cor}{\label{corollary:section-tame}} For any subinterval $I$ of $(e/(e+1),1]$, we have $\gers_{\rig}^{0} (\gerY_\rig I) \subset \pi_{2,\rig}^{-1}(\gerY_{\rig} I^{\tau})$. In fact more is true: we have
\[
\gers_{\rig}^{0} (\gerY_\rig I) =\pi_{2,\rig}^{-1}(\gerY_{\rig} I^{\tau}) \cap \pi_{1,\rig}^{-1}(\gerY_{\rig} I).
\]
\end{cor}
\begin{proof}
Since $\gers_{\rig}^{0} (\gerY_\rig I) \subset \pi_{1,\rig}^{-1}(\gerY_{\rig} I)$ and in view of Lemma \ref{lemma:pre-trace}, to show the first inclusion  we only need to show that $\gers_{\rig}^{0} (\gerY_\rig I)$ does not intersect $\pi_{2,\rig}^{-1}(\gerY_{\rig} I^{w})$. Let $Q$ be in $\gers_{\rig}^{0} (\gerY_\rig I)$ . Assumption {\bf A6} tells us that $\pi_{1,\rig}^{\prime}(Q) \in \gerY_{\rig}[0,e/(e+1))$. Therefore, $\pi_{2,\rig}(Q)=w_{\rig}\pi_{1,\rig}^{\prime}(Q) \in \gerY_{\rig}(1/(e+1),1]$. Now  $(1/(e+1),1] \cap I^{w}=\emptyset$ and the first inclusion follows. To prove the equality note that $\pi_{2,\rig}^{-1}(\gerY_{\rig} I^{\tau})$ is a subset of
\[
  \pi_{2,\rig}^{-1}(\gerY_{\rig}(1/(e+1),1]) =(\pi_{1,\rig}^\prime)^{-1}(\gerY_{\rig} [0,e/(e+1)) )=\gers^0_\rig(\gerY_{\rig} (1/(e+1),1])
\]
where for the first equality we use assumption {\bf A4} and for the second, assumption {\bf A6}. Intersecting with $\pi_{1,\rig}^{-1}(\gerY_{\rig} I)$ gives the desired result.

\end{proof}

\vspace{2mm}

 Let $\calF$ be an invertible sheaf over $X$.  Let $\gerF$
denote the induced sheaf on the formal scheme $\gerX$, and
$\gerF_\rig$ the sheaf on $\gerX_\rig$. To simplify the
 notation, we sometimes denote the sheaf $\pi^*\calF$ on $Y$ also by $\calF$.
Similarly we sometimes denote both $\pi_\rig^*\gerF_\rig$ on $\gerY_\rig$ and $\pi_{1,\rig}^*\pi_\rig^*\gerF_\rig$ on $\gerY^0_\rig$, also by $\gerF_\rig$. We try to avoid this
abbreviation when it is likely to cause confusion.  We assume that

\begin{itemize}
\item[\bf{A7}\ ]we have fixed a morphism of $\calO_Y$-modules
$\vartheta:w^*\pi^*\calF \ra \pi^*\calF$.

\end{itemize}

\begin{rem}
In the context of modular curves, the sheaf $\calF$ on $X=X(\Gamma_1(N))$ can be taken to be $\omg^{\ts k}$, where $\omg$ is the usual sheaf whose restriction to the noncuspidal locus is the push forward of the sheaf of invariant differentials of the universal elliptic curve. On $Y=X(\Gamma_1(N) \cap \Gamma_0(p))$, $w^*\pi^*\omg$ , on the noncuspidal locus, is the push forward of the sheaf of invariant differentials on the quotient of the universal elliptic curve by the distinguished subgroup of order $p$, and pulling back via this quotient morphism (raised to the power $k$) furnishes us with $\vartheta$ in this case. Similar  morphisms exist for other Shimura curves.
\end{rem}

\

\begin{defn}\label{definition:frobenius}
Let $\sigma:\gerX_\rig[0,1/(e+1)) \ra \gerX_\rig[0,e/(e+1))$ be the morphism
defined by
\[
\sigma(P)=\pi_\rig w_\rig \gers_\rig(P).
\]
Using various parts of Proposition \ref{proposition:GK}  we show  $\nu_\gerX(\sigma(P))=e\nu_\gerX(P)<e/(e+1)$.  Part (1) shows that
$\nu_\gerY(\gers_\rig(P))=\nu_\gerX(P)$. Part (4) then implies that
$\nu_\gerY(w_\rig\gers_\rig(P))=1-\nu_\gerX(P) > e/(e+1)$. Therefore by part (2) we get $\nu_\gerX(\sigma(P)) =\nu_\gerX(\pi_\rig w_\rig\gers_\rig(P))=e(1-\nu_\gerY(w_\rig\gers_\rig(P)))=e\nu_\gerX(P)$.\end{defn}

\

\begin{rem}\label{remark:analogy-frob}
In the context of modular curves, $\sigma$ corresponds to the
Frobenius morphism obtained by dividing an elliptic curve and its tame level structure by its canonical subgroup.
\end{rem}

\subsection{Overconvergence on $\gerY_\rig$.}


Let $L \subset \hat{\bar{L}}_0$ to be a completely valued extension of $L_0$ with ring of integers
${\mathcal O}$.

\begin{defn}{\label{definition:overconvergence-tame}}
Let $0 \leq r < e/(e+1)$ be in $\QQ$. The space of
$r$-overconvergent sections of $\gerF_\rig$ on $\gerY_\rig$ defined
over $L$ is
\[
S_r(\gerY_\rig,\gerF_\rig;L):=H^0(\gerY_\rig[0,r]\tsc
L,\gerF_\rig)\cong H^0(\gerX_\rig[0,r] \tsc L,\gerF_\rig)
\]
where the last identification is via $\gers_\rig^*$. An {\it
overconvergent} section of $\gerF_\rig$ on $\gerX_\rig$ defined
over $L$ is an element of $S_r(\gerY_\rig,\gerF_\rig;L)$ for some
$r>0$. The space of overconvergent sections of $\gerF_\rig$ on
$\gerX_\rig$ is denoted by $S^\dagger(\gerX_\rig,\gerF_\rig;L)$.
An element of $S_r(\gerX_\rig,\gerF_\rig;L)$ is called {\it
classical} if it is in the image of the map
\[
H^0(\gerY_\rig \tsc L,\gerF_\rig) \inc H^0(\gerY_\rig[0,r] \tsc L
,\gerF_\rig) \cong H^0(\gerX_\rig[0,r]\tsc L,\gerF_\rig)
\]
where the first arrow is restriction, and the second identification
is via $\gers_\rig^*$. The space of classical sections of
$\gerF_\rig$ over $\gerX_\rig$ defined over $L$ is denoted by ${\bf
S}(\gerX_\rig,\gerF_\rig;L)$. Note that this is possibly larger than
$H^0(\gerX_\rig \tsc L, \gerF_\rig)$. We note that if $Y$ is
projective over $\calO_0$ then by rigid analytic GAGA we know
analytification induces an isomorphism between $H^0(Y \ts L,\calF)$
and $H^0(\gerY_\rig \tsc L,\gerF_\rig)$ and hence a classical
section is indeed the analytification of an algebraic global
section.
\end{defn}

\

We define norms on these spaces making them into $p$-adic  Banach
spaces. Let $\gerZ$ be an admissible formal scheme, and $\gerZ_\rig$
its generic fibre. Let $z$ be a point of $\gerZ_\rig$. Throughout
this article the notation $\gamma_z$ stands for a choice of an
$E$-point
 of $\gerZ_\rig$ mapping onto $\{z\}$, where $E$ is a finite extension of the residue
 field of $z$. In this situation we say that the morphism
 $\gamma_z: \Sp(E) \ra \gerZ_\rig$ gives the point $z$. Any such morphism can
 be uniquely extended to a morphism $\tilde{\gamma}_z:\Spf(\calO_E)
 \ra \gerZ$, where $\calO_E$ is the ring of integers in $E$.

\

\begin{defn}{\label{definition:norm}} Let $\gerZ$ be a reduced quasi-compact
admissible formal scheme over $\calO$ and $\gerN$ an invertible
sheaf on it. Let~$z \in \gerZ_\rig$ be a point and
$\gamma_z\colon\Sp(E)\ra \gerZ_\rig$ an $E$-point giving $z$ (where
$E$ is a finite extension of the residue field of $z$). We first
define a norm $|.|_z$ on $H^0(\Sp(E),z^*\gerN_\rig)$. Denote the
formal extension of $\gamma_z$ to $\gerZ$ by
$\tilde{\gamma}_z:\Spf(\calO_{E}) \ra \gerZ$. Then
\[
H^0(\Sp(E),\gamma_z^*\gerN_\rig)=H^0(\Spf(\calO_{E}),\tilde{\gamma}_z^*\gerN)\otimes_{\calO_{E}}E
\]
and we define $|.|_z$ via identifying
$H^0(\Spf(\calO_{E}),\tilde{\gamma}_z^*\gerN)$ with $\calO_{E}$.
The definition is independent of the identification and the choice
of $\gamma_z$. Let~$\calU \subset \gerZ_\rig$ be an admissible
open, and let~$f \in H^0(\calU,\gerN_\rig)$ and $z \in \calU$. We
define
\[
|f(z)|:=|{\gamma_z}^*f|_z.
\]
We define the norm of $f$ over $\calU$ to be
$|f|_{_\calU}:=\sup\{|f(z)|: z \in \calU\}$ (possibly infinite).
\end{defn}

\

\begin{lemma}{\label{lemma:Banach}}
Assume $L$ is discretely valued. Let $\gerZ$ be a reduced quasi-compact admissible formal scheme over
$\mathcal{O}$, and $\gerN$ an invertible sheaf on $\gerZ$. If~$\
\calU$ is an affinoid subdomain of~$\gerZ_\rig$, then~$|.|_{_\calU}$
is a norm on~$H^0(\calU,\gerN_\rig)$ which makes it into a
potentially orthonormizable $L$-Banach module.
\end{lemma}

\begin{proof}
See Lemma 2.2 of \cite{Kassaei3}. The only thing left to show is the potential orthonormizability, that is, the existence of an equivalent norm on $H^0(\calU,\gerN_\rig)$ with respect to which there is an  orthonormal basis. This follow from Proposition 1 of \cite{Serre} along with remarks made before the {\it Exemple}.

\end{proof}

For simplicity, we denote $|.|_{\gerY_\rig[0,r]\tsc L}$ on
$S_r(\gerY_\rig,\gerF_\rig;L)=H^0(\gerY_\rig[0,r]\tsc
L,\gerF_\rig)$ by $|.|_r$.  It is clear that $S_r(\gerY_\rig,\gerF_\rig;L)$ is isomorphic as a normed space to $S_r(\gerY_\rig,\gerF_\rig;L_{0})\tsc L$ for any completely valued subfield $L$ of $\hat{\overline{L}}$ containing $L_{0}$.

\begin{cor}{\label{corollary:ortho}}
If $L$ is discretely valued, $S_r(\gerY_\rig,\gerF_\rig;L)$ is a potentially orthonormizable $L$-Banach module with respect to $|.|_r$.
\end{cor}

Next, we define an operator $\U$ on
$S_r(\gerY_\rig,\gerF_\rig;L)$. We start with some generalities.

\begin{defn}{\label{definition:trace}} Throughout  \S 2-4 we fix a choice of  (a ``normalization factor'') $c \in L_0$.  This will be used to define the following collection of operators.

Let $\calY$, $\calY^0$ be rigid analytic spaces over $L$, and
$\alpha_1,\alpha_2\colon\calY^0 \ra \calY$ be two finite flat
morphisms. Assume that $\calG$ is a quasi-coherent sheaf on $\calY$,
and we are given a  morphism $\ell:\alpha_2^*\calG \ra
\alpha_1^*\calG$. Let $\calV,\calW \subset \calY$ be admissible
opens such that $\alpha_1^{-1}(\calW) \subseteq
\alpha_2^{-1}(\calV)$. Then, we can define
\[
T=T_{\calV}^{\calW}: H^0(\calV,\calG) \ra H^0(\calW,\calG)
\]
via the following composite
\begin{eqnarray}{\label{diagram:trace}}
\xymatrix{ H^0(\alpha_2^{-1}(\calV),\alpha_2^*\calG) \ar[r]^{res}
&H^0(\alpha_1^{-1}(\calW),\alpha_2^*\calG)\ar[r]^{\ell}
&H^0(\alpha_1^{-1}(\calW),\alpha_1^*\calG)
\ar[d]^{(\alpha_1)_*} \\
H^0(\calV,\calG) \ar[u]^{\alpha_2^*} \ar[rr]^{T^{\calW}_{\calV}}
&&H^0(\calW,\calG)}
\end{eqnarray}
where $res$ denotes restriction, and $(\alpha_1)_*$ is the
push-forward map (i.e., the map induced by the trace map between the
structure sheaves) which is defined since $\alpha_1$ is finite and
flat.

Now, let $\calV, \calW$ be as above with the further
assumption that $\calV \subseteq \calW$. In this case, we define
an $L$-linear transformation
\[
\U_{\calV}^\calW: H^0(\calV,\calG) \ra H^0(\calV,\calG)
\]
via $\U_{\calV}^\calW=res_{\calV}^\calW \circ cT_{\calV}^{\calW}$,
where $res_{\calV}^\calW:H^0(\calW,\calG)\ra H^0(\calV,\calG)$ is
the natural restriction. When it is understood which $\calW$ we are
using, we often drop it from the notation and simply write
$\U_\calV$.
\end{defn}

\

\begin{rem}{\label{remark:after-trace}}
 Note that the above definition of $\U$ operators depends on our fixed choice of  $c$, which we have suppressed from the notation. In each specific application in \S 5,
 we will specify the choice of $c$.
For example, in the case of modular curves $c$ will be taken to be $1/p$.
\end{rem}

\

\begin{lemma}{\label{lemma:properties-of-trace}} Let $\calY, \calY_0, \alpha_1, \alpha_2$ be as above.

\begin{enumerate}

\item For $i=1,2$, let $(\calV_i,\calW_i)$ be as in Definition \ref{definition:trace}., with $\calV_1 \subset \calV_2$ and $\calW_1 \subset \calW_2$. Then
\[
res_{\calW_1}^{\calW_2}T_{\calV_2}^{\calW_2}=T_{\calV_1}^{\calW_1}res_{\calV_1}^{\calV_2}.
\]

\item Let $\calV_1 \subseteq \calV_2 \subseteq \calV_3$ be
admissible opens of $\calY$, such that $\alpha_1^{-1}(\calV_{i+1})
\subseteq \alpha_2^{-1}(\calV_{i})$ for $i=1,2$. Then we have

\vspace{2mm}

\begin{enumerate}

\item $\U^{\calV_3}_{\calV_2}=cT^{\calV_2}_{\calV_1}res^{\calV_2}_{\calV_1},$

\vspace{2mm}

\item   $res_{\calV_1}^{\calV_2}\U^{\calV_3}_{\calV_2}=\U^{\calV_2}_{\calV_1}res^{\calV_2}_{\calV_1},$

\vspace{2mm}

\item
$\U_{\calV_2}^{\calV_3}T^{\calV_2}_{\calV_1}=T^{\calV_2}_{\calV_1}   \U_{\calV_1}^{\calV_2}.$
\end{enumerate}

\end{enumerate}
\end{lemma}

\begin{proof}
Part (1) is immediate from Definition \ref{definition:trace}, and part (2) follows easily from part (1).

\end{proof}

\



Since $\pi_{2,\rig}=w_\rig\pi_{1,\rig}^\prime$, and $\pi_\rig \pi_{1,\rig}=\pi_\rig \pi_{1,\rig}^\prime$, we
can apply $(\pi^\prime_{1,\rig})^*$ to the (analytified) morphism $\vartheta$,
and get a morphism of $\calO_{\gerY^0_\rig}$-modules
\[
\vartheta:\pi_{2,\rig}^*\pi^*_\rig\gerF_\rig \ra \pi_{1,\rig}^*
\pi^*_\rig\gerF_\rig.
\]
Note that we often denote the sheaf on the right hand side simply by $\gerF_\rig$.

\

\begin{defn}{\label{definition:U}}  Let $0 \leq r < 1/(e+1)$. Following the notation of Definition \ref{definition:trace} we let $\calY=\gerY_\rig$,
$\calY^0=\gerY^0_\rig$, $\alpha_i=\pi_{i,\rig}$ for $i=1,2$,
$\calG=\gerF_\rig$, $\calW=\gerY_\rig[0,er]$,
$\calV=\gerY_\rig[0,r]$, and $\ell=\vartheta$. In view of Lemma \ref{lemma:pre-trace}     we have
\[
\pi_{1,\rig}^{-1}(\calW) \subseteq
\pi_{2,\rig}^{-1}(\calV).
\]
Hence  Definition \ref{definition:trace} gives us  $T=T^{\calW}_{\calV}=T^{\gerY_\rig[0,er]}_{\gerY_\rig[0,r]}$, and also (after tensoring with $L$)
\[
\U_r:=\U^{\gerY_\rig[0,er]}_{\gerY_\rig[0,r]}:S_r(\gerY_\rig,\gerF_\rig;L)
\ra S_r(\gerY_\rig,\gerF_\rig;L).
\]
It is easy to see (for example using part (1) of Lemma \ref{lemma:properties-of-trace}) that  if $0 \leq r \leq r^\prime <  1/(e+1)$ then $\U_r$ induces an operator on $S_{r^\prime}(\gerY_\rig,\gerF_\rig;L)$ which equals $S_{r^\prime}$. Therefore, we obtain an operator
\[
\U:S^\dagger(\gerY_\rig,\gerF_\rig;L) \ra
S^\dagger(\gerY_\rig,\gerF_\rig;L).
\]
We sometimes denote $\U_r$ by $\U$ as well.

\end{defn}

\

\begin{prop} {\label{proposition:compcont-tame}} Assume $L$ is discretely valued. The operator $\U_r$ is a completely continuous operator
for any  $0< r < 1/(e+1)$.
\end{prop}

\begin{proof} Since by definition $\U_r$ is the composite of a continuous
operator $T$ with restriction of sections from $\gerY_\rig[0,er]$ to
$\gerY_\rig[0,r]$ it is enough to show that this restriction is a
morphism of Banach spaces with compact closure. Since $r<1/(e+1)$
using the section $\gers_\rig$ it is enough to show that restriction
from $H^0(\gerX_\rig[0,er],\gerF_\rig)$ to
$H^0(\gerX_\rig[0,r],\gerF_\rig)$ has compact closure. This can be
shown, for example, by using Propositions 2.4.1 and 2.3.2 of
\cite{KisinLai}.   In the proof of Lemma \ref{lemma:affinoid} we
explained that (for an appropriate choice of a Cartier divisor)
$\gerX_\rig(p^{-r})$ as defined in \S 2.3 of \cite{KisinLai} is the
same as $\gerX_\rig[0,r]$ in our notation. Now we can apply
Proposition 2.4.1 of \cite{KisinLai} noting that the Banach space
structure we have defined on $H^0(\gerX_\rig[0,r],\gerF_\rig)$
equals its canonical Banach space structure.  \end{proof}

\

\subsection{Set-up: the case of ``higher levels''}

In this section we axiomatize the situation common to Shimura curves
with level structures containing ``powers of $p$". Fix a positive
integer $m$.  Assume $\Xm$ is a curve over $\calO_0$. Let $\gerXm$
denote the formal completion of $\Xm$ along its special fibre.
Denote the rigid analytic generic fibre of this quasi-compact
admissible formal scheme by $\gerXm_\rig$. Let $\gerXmo_\rig$ be a
rigid analytic curve over $L_{0}$  (the notation is chosen for
uniformity and is not meant to suggest that $\gerXmo_\rig$  is the
rigid analytic fibre of a canonically chosen formal scheme). Assume
we have two morphisms
$\lambda_{1,\rig},\lambda_{2,\rig}\colon\gerXmo_\rig \ra
\gerXm_\rig$ defined over $L_{0}$. Assume further that
\begin{itemize}

\item[\bf{H1}\ ] $\Xm \ts L_0$ is smooth;

\item[\bf{H2}\ ] $\lambda_{1,\rig}$, $\lambda_{2,\rig}$ are finite flat
morphisms;

\item[\bf{H3}\ ] for $1 \leq i \leq m$, there are morphisms
$\phi_i:\Xm \ra Y$ and $\eta_{i,\rig}:\gerXmo_\rig \ra \gerY^0_\rig$ fitting in the following commutative diagrams:

\begin{itemize}

\item[\bf{H3.1}\ ] for each $ 1 \leq i \leq m$, and $j=1,2$ the
diagrams

\begin{eqnarray}\label{diagram:higherlevel}
\xymatrix{ \gerXmo_\rig\ar[r]^{\lambda_{j,\rig}} \ar[d]_{\eta_{i,\rig}} & \gerXm_\rig
\ar[d]^{\phi_{i,\rig}}  \\
\gerY_\rig^0 \ar[r]_{\pi_{j,\rig}} &\gerY_\rig}
\end{eqnarray}
are commutative. Furthermore, for $i=j=1$
the diagram obtained is a {\underline{product}} diagram.

\item[\bf{H3.2}] For each $ 1 \leq i \leq m-1$, the diagram

\begin{eqnarray}\label{diagram:higherlevel2}
\xymatrix{ \gerXmo\ar[r]^{\lambda_{1,\rig}} \ar[d]_{\eta_{i+1,\rig}} & \gerXm_\rig
\ar[d]^{\phi_{i,\rig}}  \\
\gerY_\rig^0 \ar[r]_{\delta_\rig\pi_{2,\rig}} &\gerY_\rig}
\end{eqnarray}
is commutative.

\end{itemize}

\end{itemize}

\begin{rem} Let us continue to hint at what the above curves and
morphisms will signify in the context of Shimura curves by
presenting them in the case of modular curves. In that context,
$\Xm$ can be taken to be $X(\Gamma_1(N) \cap \Gamma_1(p^m))$, whose
noncuspidal points classify $(E,i,P)$ where $(E,i)$ is an elliptic
curve with  $\Gamma_1(N)$-level  structure and $P$ is a point of
exact order $p^m$ on $E$ (in the sense of Drinfeld). Then
$\gerXmo_\rig$ will be the analytification of the curve over
$L_0=\QQ_p$ whose noncuspidal points classify $(E,i,P,D)$ where one
adds $D$, a subgroup of order $p$, to the above data assuming it
intersects the subgroup generated by $P$ trivially. The morphisms
$\lambda_{1,\rig}$ and $\lambda_{2,\rig}$ are the analytification of
morphisms which forget $D$ and quotient out by $D$, respectively.
The morphism $\phi_j$ can be taken to send $(E,i,P)$ to $(E/\langle
p^{m+1-j}P\rangle, \bar{i}, \langle p^{m-j}P \rangle)$. The morphism
$\eta_{i,\rig}$ can be defined in the same way by further enclosing
the image of $D$ in the quotient. All the required properties can be
easily checked to hold using Yoneda's lemma.
\end{rem}

We will again denote by $\calF$ the sheaf $\phi_1^*\calF$ on $\Xm$,
and by $\gerF_\rig$ the corresponding sheaf on $\gerXm_\rig$. We
will also denote by $\gerF_\rig$  the sheaf
$\eta_{1,\rig}^*\gerF_\rig$ on $\gerXmo_\rig$.  Using {\bf H3.1} with $i=j=1$, it can also
be described as $\lambda_{1,\rig}^*\gerF_\rig=\lambda_{1,\rig}^*\phi_{1,\rig}^*\gerF_\rig$, notations that we may sometimes use. Applying $\eta_{1,\rig}^*$ to the morphism of
$\calO_{\gerY^0_\rig}$-modules $\vartheta:\pi_{2,\rig}^*\gerF_\rig
\ra \pi_{1,\rig}^*\gerF_\rig$, we obtain a morphism of
$\calO_{\gerXmo_\rig}$-modules
\[
\lambda_{2,\rig}^*\gerF_\rig \ra \lambda_{1,\rig}^*\gerF_\rig
\]
which we still denote by $\vartheta$.

\

We will use the following lemma later.

\begin{lemma}\label{lemma:redundant} Let $Q \in \gerXm_\rig$, and $1\leq i \leq m-1$. If
$(\phi_{i,\rig}(Q))^w$ is canonical, then $\phi_{i+1,\rig}(Q)$ is
anti-canonical.
\end{lemma}

\begin{proof} We have $\nu_{\gerY}(\phi_{i,\rig}(Q))=1-\nu_{\gerY}((\phi_{i,\rig}(Q))^w)>1/(e+1)$ by Proposition \ref{proposition:GK}. Let $R \in \gerXmo_\rig$ be a point such that
$\lambda_{1,\rig}(R)=Q$. Then
\[
\delta_\rig w_\rig\pi_{1,\rig}^\prime
\eta_{i+1,\rig}(R)=\delta_\rig\pi_{2,\rig}\eta_{i+1,\rig}(R)=\phi_{i,\rig}\lambda_{1,\rig}(R)=\phi_{i,\rig}(Q)
\]
by assumption {\bf H3.2}. Proposition \ref{proposition:GK} now implies
$\nu_\gerY(\pi_{1,\rig}^\prime
\eta_{i+1,\rig}(R))=1-\nu_\gerY(\phi_{i,\rig}(Q))<e/(e+1)$.
Therefore $\pi_{1,\rig}^\prime \eta_{i+1,\rig}(R)$ is canonical and
hence, by Assumption {\bf A5}, we know that $\pi_{1,\rig}
\eta_{i+1,\rig}(R)$ is anti-canonical. But Assumption {\bf H3.1}
tells us that
\[
\pi_{1,\rig}
(\eta_{i+1,\rig}(R))=\phi_{i+1,\rig}(\lambda_{1,\rig}(R))=\phi_{i+1,\rig}(Q).
\]
We are done.
\end{proof}

\

\subsection{Overconvergence over $\gerXm_\rig$} In this subsection we define the space of overconvergent
sections of $\gerF_\rig$ on $\gerXm_\rig$. Recall that  $L \subset \hat{\bar{L}}_0$ is a completely valued extension of $L_0$ with ring of integers ${\mathcal O}$.
\

\begin{defn}{\label{definition:overconvergence-higher}}
Let $0 \leq r < e^{2-m}/(e+1)$ be in $\QQ$. The space of
$r$-overconvergent sections of $\gerF_\rig$ on $\gerX(m)_\rig$
defined over $L$ is
\[
S_r(\gerXm_\rig,\gerF_\rig;L):=H^0{\large
(}\phi_{m,\rig}^{-1}(\gerY_\rig[0,re^{m-1}] \tsc
L),\gerF_\rig{\large )}.
\]
An overconvergent section of $\gerF_\rig$ on $\gerXm_\rig$ is an
element of $S_r(\gerXm_\rig,\gerF_\rig;L)$ for some $r>0$. The
space of overconvergent sections of $\gerF_\rig$ on $\gerXm_\rig$
is denoted by $S^\dagger(\gerXm_\rig,\gerF_\rig;L)$. An
overconvergent section of $\gerF_\rig$ is called {\it classical}
if it is in the image of the restriction map
\[
H^0(\gerXm_\rig \tsc L,\gerF_\rig) \inc H^0{\large
(}\phi_{m,\rig}^{-1}(\gerY_\rig[0,re^{m-1}] \tsc
L),\gerF_\rig{\large )}.
\]
We denote the space of these classical sections by ${\bf
S}(\gerXm_\rig,\gerF_\rig;L)$.  If $X(m)$ is projective over $\calO_0$ then by rigid analytic GAGA  analytification induces an isomorphism between $H^0(X(m) \ts L,\calF)$ and $H^0(\gerXm_\rig \tsc L,\gerF_\rig)$
and hence a classical section is indeed the analytification of an algebraic global section.
\end{defn}

\

By Definition \ref{definition:norm} we have a norm
$|.|_{_{\phi_{m,\rig}^{-1}(\gerY_\rig[0,re^{m-1}])}}$ on
$S_r(\gerXm_\rig,\gerF_\rig;L)$ which we abbreviate by $|.|_r$. As
in Corollary \ref{corollary:ortho} we have the following.
\begin{prop}
If $L$ is discretely valued, the space $S_r^m(\gerXm_\rig,\gerF_\rig;L)$ is a potentially orthonormizable $L$-Banach module with respect to $|.|_r$.
\end{prop}

\

Now we prepare for the definition of the $\U$ operator. Let $0 \leq r < e^{1-m}/(e+1)$. Let
$\calW=\phi_{m,\rig}^{-1}(\gerY_\rig[0,re^{m}])$, and
$\calV=\phi_{m,\rig}^{-1}(\gerY_\rig[0,re^{m-1}])$. Since $0\leq re^{m}<e/(e+1)$, we can apply part (1) of Lemma \ref{lemma:pre-trace} to deduce that
\[
\pi_{1,\rig}^{-1}(\gerY_\rig[0,re^{m}])\subseteq
\pi_{2,\rig}^{-1}(\gerY_\rig[0,re^{m-1}]).
\]
Applying $\eta_{m,\rig}^{-1}$ to both sides, and noting that by
Assumption {\bf H3.2} we have  $\phi_{m,\rig}\lambda_{i,\rig}=\pi_{i,\rig}\eta_{m,\rig}$ for $i=1,2$ we get
\begin{eqnarray}\label{equation:pre-U-higher}
\lambda_{1,\rig}^{-1}(\calW)\subseteq
\lambda_{2,\rig}^{-1}(\calV).
\end{eqnarray}

\begin{defn} {\label{definition:U-higher}}Assume $0 \leq r < e^{1-m}/(e+1)$. Following the notation of Definition \ref{definition:trace} we set $\calY=\gerXm_\rig$,
$\calY^0=\gerXmo_\rig$, $\alpha_i=\lambda_{i,\rig}$ for $i=1,2$,
$\calF=\gerF_\rig$,
$\calW=\phi_{m,\rig}^{-1}(\gerY_\rig[0,re^{m}])$,
$\calV=\phi_{m,\rig}^{-1}(\gerY_\rig[0,re^{m-1}])$, and
$\ell=\vartheta$. By Equation \ref{equation:pre-U-higher} we can define
\[
T=T_{\calV}^{\calW}=T_{\phi_{m,\rig}^{-1}(\gerY_\rig[0,re^{m-1}])}^{\phi_{m,\rig}^{-1}(\gerY_\rig[0,re^{m}])}
\]
and in turn (after tensoring with $L$)
\[
\U_r:=\U_{\calV}^{\calW}: S_r(\gerXm_\rig,\gerF_\rig;L) \ra
S_r(\gerXm_\rig,\gerF_\rig;L)
\]
As in Definition \ref{definition:U} we get an operator
\[
\U_r: S^\dagger(\gerXm_\rig,\gerF_\rig;L) \ra
S^\dagger(\gerXm_\rig,\gerF_\rig;L)
\]
We often denote $\U_r$ simply by $\U$.
\end{defn}

\
 As in Proposition \ref{proposition:compcont-tame} we can prove the following.

\begin{prop}{\label{proposition:compcont-higher}}  Assume $r>0$ and $L$ is discretely valued. Then $\U_r$ is a completely continuous
operator of $S_r(\gerXm_\rig,\gerF_\rig;L)$.
\end{prop}

\

\section{Analytic Continuation}

\subsection{Generalities on Analytic continuation.}
The results of this section and their applications in the subsequent sections \S \ref{subsection:ancont}, \S \ref{subsection:ancont-higher} are based on Buzzard's method of analytic continuation in \cite{Buzzard}. We must note that we avoid requiring any connectedness properties for the regions involved in this process of analytic continuation as in the examples provided by Shimura curves there are no cusps: over modular curves one can show that two sections agree on a connected region by showing that their $q$-expansions agree.

Let $L$ be a completely valued subfield of  ${\bar{L}}_0$.
Let $\calY$, $\calY^0$ be rigid analytic spaces over $L$, and
$\alpha_1,\alpha_2\colon\calY^0 \ra \calY$ be two finite flat
morphisms. Assume that $\calG$ is a quasi-coherent sheaf on $\calY$,
and we are given a morphism $\ell:\alpha_2^*\calG \ra
\alpha_1^*\calG$.

Let $\calU_0$ be an admissible open of $\calY$. Assume we are given
an admissible open $\calU$ of $\calY$, and an admissible covering
\[
\calU_0 \subseteq \calU_1 \subseteq \cdots \subseteq \calU_n
\subseteq \cdots
\]
of $\calU$ such that $\alpha_1^{-1}(\calU_{n+1}) \subseteq
\alpha_2^{-1}(\calU_{n})$ for all $n \geq 0$.

For each $n \geq 0$, let
\[
\U_n:H^0(\calU_n,\calG) \ra H^0(\calU_n,\calG)
\]
be the operator denoted by $\U_{\calU_{n}}^{\calU_{n+1}}$ in
Definition \ref{definition:trace}.  From the above
assumptions it follows that $\alpha_1^{-1}(\calU) \subseteq
\alpha_2^{-1}(\calU)$, and hence we have an operator
\[
\U_\infty:H^0(\calU,\calG)\ra H^0(\calU,\calG)
\]
which is the operator $T_{\calU}^\calU=\U_{\calU}^{\calU}$ in the notation of
Definition \ref{definition:trace}.

\begin{prop}{\label{proposition:analytic-general}}
Let the notation be as above. Let $f_0\in H^0(\calU_0,\calG)$.
\begin{enumerate}
\item Let $A(x) \in L[X]$ be a polynomial such that $A(0) \neq 0$.
If $A(\U_0)f_0$ can be extended to a section $F$ of $\calG$ on
$\calU$, then so can $f_0$. Denote this extension of $f_0$ by $f
\in H^0(\calU, \calG)$. We have $A(\U_\infty)f=F$.

\item (Special case of the first part). If $f_0$ is a generalized eigensection for $\U_0$ of generalized eigenvalue $a \neq 0$, then $f_0$ can be extended to a section $f\in
H^0(\calU,\calG)$ which is a generalized eigensection for $\U_\infty$ with generalized eigenvalue $a$.
\end{enumerate}

\end{prop}

\begin{proof}
Let $A(x)=xA_0(x)-a$, where $0 \neq a \in L$.  For $n \geq 1$ define
$f_n \in H^0(\calU_n,\calG)$, recursively, via
\begin{eqnarray}{\label{equation:recursive}}
f_{n}:=a^{-1}
cT^{\calU_n}_{\calU_{n-1}}A_0(\U_{n-1})f_{n-1}-a^{-1}F_{|_{\calU_n}}.
\end{eqnarray}
To prove the statement, it is enough to prove
\begin{eqnarray}
&{\mathbf{P}(n)}\ \ \ \ \ &A(\U_{n-1})f_{n-1}=F_{|_{\calU_{n-1}}}\nonumber\\
&{\mathbf{Q}(n)}\ \ \ \ \ &{f_{n}}_{|_{\calU_{n-1}}}=f_{n-1}
\end{eqnarray}
for all $n \geq 1$. The reason is that knowing $\mathbf{Q}(n)$ for
all $n$ would give us a section $f$ of $\calG$ on $\calU$ such
that $f_{|_{\calU_n}}=f_n$ for all $n\geq 0$, and knowing
$\mathbf{P}(n)$ for all $n$ would imply that for all $n \geq 0$ we
have
\[
(A(\U_\infty)f)_{|_{\calU_n}}=A(\U_n)(f_{|_{\calU_n}})=A(\U_n)(f_n)=F_{|_{\calU_n}}
\]
where the first equality holds by part (2)(b)
of Lemma \ref{lemma:properties-of-trace}.

On the other hand, using Equation \ref{equation:recursive}, it is easy to see that $\mathbf{P}(n)$ implies
 $\mathbf{Q}(n)$ for all $n\geq 1$. Hence it suffices to prove $\mathbf{P}(n)$ for
all $n\geq 1$. We do this by induction. By assumption
$\mathbf{P}(1)$ holds. Assume $\mathbf{P}(n)$ holds. We have
\begin{eqnarray}
A(\U_{n})f_{n}=&\nonumber\\
A(\U_{n}){\large(}a^{-1}cT^{\calU_{n}}_{\calU_{n-1}}A_0(\U_{n-1})f_{n-1}-a^{-1}F_{|_{\calU_n}}{\large)}=&\nonumber\\
a^{-1}cT^{\calU_n}_{\calU_{n-1}}A(\U_{n-1})A_0(\U_{n-1})f_{n-1}-a^{-1}A(\U_{n})F_{|_{\calU_n}}=&\nonumber\\
a^{-1}cT^{\calU_n}_{\calU_{n-1}}A_0(\U_{n-1})A(\U_{n-1})f_{n-1}-a^{-1}A(\U_{n})F_{|_{\calU_n}}=&\nonumber\\
a^{-1}cT^{\calU_n}_{\calU_{n-1}}A_0(\U_{n-1})(F_{|_{\calU_{n-1}}})-a^{-1}A(\U_{n})F_{|_{\calU_n}}=&\nonumber\\
a^{-1}A_0(\U_n)cT^{\calU_n}_{\calU_{n-1}}(F_{|_{\calU_{n-1}}})-a^{-1}A(\U_{n})F_{|_{\calU_n}}=&\nonumber\\
 a^{-1}A_0(\U_n) \U_{n}(F_{|_{\calU_n}})-a^{-1}A(\U_{n})F_{|_{\calU_{n}}}=&F_{|_{\calU_{n}}},\nonumber\\
\end{eqnarray}
where for the second and fifth equalities we use part (2)(c) of Lemma
\ref{lemma:properties-of-trace}, and for the sixth equality we use part (2)(a) of the same lemma. The second part of the proposition follows from the first part by setting
$A(x)=(x-a)^N$ for some $N \geq 1$, and $F=0$.
\end{proof}

\subsection{Analytic continuation on $\gerY_\rig$.}{\label{subsection:ancont}}

\begin{prop}{\label{proposition:analytic-tame}} Let $0<r<1/(e+1)$ be in $\QQ$. Let $f\in S_r(\gerY_\rig,\gerF_\rig;L)=
H^0(\gerY_\rig[0,r]\tsc L,\gerF_\rig)$.
\begin{enumerate}
\item Let $A(x) \in L[X]$ be a polynomial such that $A(0) \neq 0$.
If $A(\U)f$ can be extended to a section $F$ of $\gerF_\rig$ on
$\gerY_\rig[0,1)\tsc L$, then so can $f$.

\item (Special case of the first part). If $f$ is a generalized
$\U$-eigensection with  a non-zero generalized eigenvalue, then $f$ can be extended to $\gerY_\rig[0,1)\tsc L$.
\end{enumerate}

\end{prop}

\begin{proof}   For simplicity in the notation we assume, without loss of generality, that $L=L_0$. Choose $N \geq 1 $ such that
$1/(e+1)\leq e^Nr < e/(e+1)$. For $0\leq n \leq N$ define
$\calW_n=\gerY_\rig[0,e^{n-N}/(e+1)]$. We have
\[
\calW_0 \subseteq  \calW_1 \subseteq \cdots \subseteq \calW_N=\gerY_\rig[0,1/(e+1)]
\]
By part (1) of Lemma \ref{lemma:pre-trace}, for all $0\leq n \leq N-1$, we have
\[
\pi_{1,\rig}^{-1}(\calW_{n+1})\subseteq
\pi_{2,\rig}^{-1}(\calW_{n}).
\]
For  $n \geq 0$ define
$\calV_n:=\gerY_\rig[0,1-1/e^{n-1}(e+1)]$. We have
\[
\calW_0 \subseteq  \calW_1 \subseteq \cdots \subseteq \calW_N=\gerY_\rig[0,1/(e+1)]=\calV_0 \subseteq \calV_1 \subseteq \calV_2 \cdots.
\]
providing an admissible covering of $\gerY_\rig[0,1)$ by affinoids (See Lemma \ref{lemma:affinoid}). Breaking up $\calV_{n+1}$ as
\[
\gerY_\rig[0,e/(e+1)) \cup \gerY_\rig[e/(e+1),e/(e+1)] \cup \gerY_\rig(e/(e+1),1-1/e^{n}(e+1)]
\]
(with the last term being empty when $n=0$) and applying all the three parts of Lemma \ref{lemma:pre-trace} we see that
\[
\pi_{1,\rig}^{-1}(\calV_{n+1})\subseteq
\pi_{2,\rig}^{-1}(\calV_n)
\]
for all $n \geq 0$. Now define $\calU_i=\calW_i$ for $i=0,\cdots,N$, and $\calU_i=\calV_{i-N}$ for $i>N$.  We have shown that for all $n \geq 0$
 \[
\pi_{1,\rig}^{-1}(\calU_{n+1})\subseteq
\pi_{2,\rig}^{-1}(\calU_n)
\]
and hence we can apply Proposition \ref{proposition:analytic-general} with the admissible covering $\{\calU_i\}$ of $\gerY_\rig[0,1)$ and $\alpha_j=\pi_{j,\rig}$ for $j=1,2$ to deduce that $f$ can be extended to a section of $\gerF_\rig$ on  $\gerY_\rig[0,1)$. Denote this extension by $f$ again. We have $A(\U)f=F$, where
\[
\U=\U^{\gerY_\rig[0,1)}_{\gerY_\rig[0,1)}
\]
in the notation of Definition \ref{definition:trace}. We are done.

\end{proof}

\begin{rem}{\label{remark:after-analytic-tame}} Note that  we have in fact shown that   $\pi_{1,\rig}^{-1} (\gerY_{\rig}[0,1)) \subseteq \pi_{2,\rig}^{-1}(\gerY_{\rig}[0,1))$ (by letting $n$ go to $\infty$), a fact implicit in the notation $\U^{\gerY_\rig[0,1)}_{\gerY_\rig[0,1)}$.

\end{rem}

\subsection{Analytic Continuation On $\gerXm_\rig$.}{\label{subsection:ancont-higher}}

\begin{lemma}{\label{lemma:inclusion}} We have
\[
\phi_{m,\rig}^{-1}(\gerY_\rig[0,1)) \subseteq
\phi_{m-1,\rig}^{-1}(\gerY_\rig[0,1)) \subseteq \cdots \subseteq
\phi_{1,\rig}^{-1}(\gerY_\rig[0,1)).
\]
\end{lemma}

\begin{proof}Let~$1\!\leq i\!\leq\!m-1$ and   $Q\!\in\!\phi_{i+1,\rig}^{-1}(\gerY_\rig[0,1))$. Assume $Q=\lambda_{1,\rig}(R)$. By assumption {\bf H3.1 }
\[
Q_1:=\pi_{1,\rig}(\eta_{i+1,\rig}(R))=\phi_{i+1,\rig}(\lambda_{1,\rig}(R))\in \gerY_\rig[0,1).
\]
We assume $\nu_\gerY(\phi_{i,\rig}(Q))=1$ and draw a contradiction. By assumption {\bf H3.2}
\[
\phi_{i,\rig}(Q)=\phi_{i,\rig}(\lambda_{1,\rig}(R))=\delta_\rig (\pi_{2,\rig}( \eta_{i+1,\rig}(R)))=\delta_\rig(w_\rig (\pi_{1,\rig}^\prime (\eta_{i+1,\rig}(R))))
\]
Hence using  parts (4) and (5) of Proposition \ref{proposition:GK} we get $\nu_\gerY( \pi_{1,\rig}^\prime(\eta_{i+1,\rig}(R)))=0$. Let $P=\eta_{i+1,\rig}(R)$. We have shown that $Q_2:=\pi_{1,\rig}^\prime(P)$ has $\nu_\gerY=0$ and hence is canonical. Therefore assumption {\bf A5} tells us that $Q_1=\pi_{1,\rig}(P)$ is anti-canonical.  Since $Q_1$ and $Q_2$ have the same image under $\pi_\rig$,  parts (1) and (2) of Proposition \ref{proposition:GK} imply that $\nu_\gerY(Q_1)=1$ which is a contradiction.

 \end{proof}

\begin{lemma}{\label{lemma:pre-analytic-higherlevel}} For each $1\leq i \leq m-1$, we have
\[
\lambda_{1,\rig}^{-1}(\phi_{i,\rig}^{-1}(\gerY_\rig[0,1))   )
=
\lambda_{2,\rig}^{-1}(\phi_{i+1,\rig}^{-1}(\gerY_\rig[0,1)))
\]

\end{lemma}
\begin{proof}
By Assumptions {\bf H3.1, \bf H3.2} we have
$\phi_{i,\rig}\lambda_{1,\rig}=\delta_\rig\pi_{2,\rig}\eta_{i+1,\rig}=\delta_\rig\phi_{i+1,\rig}\lambda_{2,\rig}$, and the result follows in view of part (5) of Proposition \ref{proposition:GK}.
\end{proof}

\begin{prop}{\label{proposition:analytic-higher}} Let $0<r<e^{2-m}/(e+1)$ be in $\QQ$.
Let $f\in S_r(\gerXm_\rig,\gerF_\rig;L)$.
\begin{enumerate}
\item Let $A(x) \in L[X]$ be a polynomial such that $A(0) \neq 0$.
If $A(\U)f$ can be extended to a section of $\gerF_\rig$ on
$\phi_{1,\rig}^{-1}(\gerY_\rig[0,1))\tsc L$, then so can $f$.

\item (Special case of the first part). If $f$ is a generalized
$\U$-eigensection with non-zero generalized eigenvalue, then $f$ can be extended to $\phi_{1,\rig}^{-1}(\gerY_\rig[0,1))\tsc L$.
\end{enumerate}

\end{prop}

\begin{proof} Without loss of generality, we can assume $L=L_0$. By definition $f$ is a section of $\gerF_\rig$ over $\phi_{m,\rig}^{-1}(\gerY_\rig[0,re^{m-1}])$. First we show that under the assumptions  $f$ extends to
$\phi_{m,\rig}^{-1}(\gerY_\rig[0,1))$. The proof of this part is exactly like the proof of Proposition \ref{proposition:analytic-tame}. One needs only apply $\phi_{m,\rig}^{-1}$ and use appropriate commutative diagrams. We will be using the notation $\calU_{i}$ from that proof.
Define $\calU^{\prime}_{i}=\phi_{m,\rig}^{-1}(\calU_i)$. Then $\{\calU^{\prime}_{i}\}$ is an admissible covering for $\phi_{m,\rig}^{-1}(\gerY_{\rig}[0,1))$. Since for $j=1,2$ we have by assumption  {\bf H3.1}
\[
\lambda_{j,\rig}^{-1} \phi_{m,\rig}^{-1}=\eta_{m,\rig}^{-1}\pi_{j,\rig}^{-1}
\]
applying $\eta_{m,\rig}^{-1}$  to the inclusion $\pi_{1,\rig}^{-1}(\calU_{i+1}) \subseteq \pi_{2,\rig}^{-1}(\calU_{i})$ for $i \geq 0$ (verified in the proof of Proposition \ref{proposition:analytic-tame}) shows that assumptions of Proposition \ref{proposition:analytic-general} are valid for the covering $\{\calU^{\prime}_i\}$, with $\alpha_{j}=\lambda_{j,\rig}$ for $j=1,2$.
 Therefore, $f$ can be extended to a section of $\gerF_{\rig}$ on $\phi_{m,\rig}^{-1}(\gerY_{\rig}[0,1))$ (which we continue to denote by $f$). Furthermore, we have $A(\U)f=F$ with
\[
\U=\U_{\phi_{m,\rig}^{-1}(\gerY_{\rig}[0,1))}^{\phi_{m,\rig}^{-1}(\gerY_{\rig}[0,1))}
\]
as in Definition \ref{definition:trace}. The final step is to extend $f$ from to $\phi_{m,\rig}^{-1}(\gerY_\rig[0,1))$ to
$\phi_{1,\rig}^{-1}(\gerY_\rig[0,1))$ (note that this statement
makes sense by Lemma \ref{lemma:inclusion}). Let us denote
$\phi_{n,\rig}^{-1}(\gerY_\rig[0,1))$ by $\calU^{\prime\prime}_{m-n}$ for $1
\leq n \leq m$.  Then, by Lemma
\ref{lemma:inclusion}, we have an admissible covering of
$\phi_{1,\rig}^{-1}(\gerY_\rig[0,1))$ given by
\[
\calU^{\prime\prime}_0 \subseteq \calU^{\prime\prime}_1 \subseteq \cdots \subseteq
\calU^{\prime\prime}_{m-1}=\phi_{1,\rig}^{-1}(\gerY_\rig[0,1))=\phi_{1,\rig}^{-1}(\gerY_\rig[0,1))=\cdots
\]
Lemma \ref{lemma:pre-analytic-higherlevel} and Remark \ref{remark:after-analytic-tame} (in conjunction with an application of $\eta_{1,\rig}^{-1}$ in the same way as above) allow us to apply Proposition \ref{proposition:analytic-general} to conclude that $f$ can be extended to a section (again denoted $f$) of $\gerF_{\rig}$ on $\phi_{1,\rig}^{-1}(\gerY_\rig[0,1))$. Proposition \ref{proposition:analytic-general} also tells us that we have
$A(\U)f=F$, where $\U$ stands for $\U_{\phi_{1,\rig}^{-1}(\gerY_\rig[0,1))}^{\phi_{1,\rig}^{-1}(\gerY_\rig[0,1))}$ here.
\end{proof}

\section{Classicality}

\subsection{The weights of $\gerF_\rig$.} We define the weights of $\calF$ on $X$.
\

\begin{defn}\label{definition:omg}
We say that $\calF$ has  a weight $k \in \ZZ$ if $\calF=\omg^{\ts k}$ for
some invertible sheaf $\omg$ on $X$, such that
$\vartheta=(\pr^*)^{\ts k}$ (see assumption {\bf A7}) where $\pr^*:w^*\pi^*\omg \ra
\pi^*\omg$ is a morphism of invertible sheaves such that
\begin{itemize}
\item[\bf A8] the morphism of $\calO_{X \ts \kappa}$-modules
\[ s^* (\pr^* \ts \kappa):s^*(w \ts \kappa)^*(\pi \ts
\kappa)^*(\omg\ts \kappa) \ra \omg \ts \kappa
\]
is the {\underline{zero}} morphism.
\end{itemize}
\noindent In what follows, for simplicity, we often follow the same notational convention involving $\calF$: denoting still by $\omg$, the sheaves $\pi^* \omg$ on $Y$, $\pi_{1,\rig}^* \pi_\rig^* \omg$ on $\gerY_\rig^0$, etc.


We define the weights of $\gerF$ and $\gerF_\rig$ to be the same as those of $\calF$. If $k$ is a weight for all these sheaves, then any section of any of these sheaves over any open set in any of the curves introduced so far is also said to have a weight $k$. We will denote the induced sheaf $\omg_\rig$ on $\gerX_\rig$ simply by $\omg$. Throughout the rest of this section we will fix a weight $k$ and a morphism $\pr^*$ satisfying Assumption {\bf A8} for $\calF$.

\end{defn}


\ For the notation in this passage, we refer the reader to the
paragraph before Definition \ref{definition:norm}. Fix a point $P$
of $\gerX_\rig[0,1/(e+1))$ for this discussion, and choose
$\gamma_P:\Sp(E)\ra \gerX_\rig[0,1/(e+1))$ giving $P$ . Then
$\gamma_{\gers_\rig(P)}:=\gers_\rig\gamma_{P}:\Sp(E)\ra
\gerY_{\rig}$ gives the point $\gers_{\rig}(P)$. Also
$\gamma_{\sigma(P)}:=\pi_{\rig}w_{\rig}\gers_{\rig}\gamma_{P}:\Sp(E)
\ra \gerX_{\rig}$ gives $\sigma(P)$ (see Definition
\ref{definition:frobenius}).  Let  $\tilde{\gamma}_{P}$,
$\tilde{\gamma}_{\gers_{\rig}(P)}$, and $\tilde{\gamma}_{\sigma(P)}$
denote the formal extensions of these maps. Then by uniqueness of
formal extensions we have
\begin{eqnarray}
\tilde{\gamma}_{P}&=&\pi\tilde{\gamma}_{\gers_{\rig}(P)}\\
\tilde{\gamma}_{\sigma(P)}&=&\pi w\tilde{\gamma}_{\gers_{\rig}(P)}.
\end{eqnarray}
The morphism $\pr^*$ given in Definition
\ref{definition:omg} induces a morphism of rigid analytic sheaves
on $\gerY_\rig$
\[
\pr^*:w_\rig^*\pi_\rig^*\omg \ra \pi_\rig^*\omg.
\]
Applying $ \gamma^{*}_{\gers_{\rig}(P)}$ to this morphism gives an
$E$-linear morphism
\begin{eqnarray}\label{equation:pr_x-rigid-sigma}
\ \ \ \pr_{\gamma_P}^*:\gamma_{\sigma(P)}^* \omg \ra
\gamma_P^*\omg.\nonumber
\end{eqnarray}
We can give an integral model for the morphism of sheaves $\pr_{\gamma_P}^*$ by  applying $\tilde{\gamma}_{\gers_\rig(P)}^*$ to (the formal
completion of) $\pr^*$ (in Definition \ref{definition:omg}). In view of   Equations 4.1 and 4.2, this gives an $\calO_E$-linear morphism
\[
\pr_{\tilde{\gamma}_P}^*: \tilde{\gamma}_{\sigma(P)} ^*  \omg \ra
\tilde{\gamma}_P^* \omg
\]
which is the restriction of $\pr_{\gamma_P}^*$ to
$\tilde{\gamma}_{\sigma(P)}^*  \omg \subset \gamma_{\sigma(P)}^*
\omg$ (where the $\omg$ on the left is over $\gerX$ and the $\omg$ on the right is over $\gerX_{\rig}$). Let $\mu$ be an element of $\calO_E$ such that $|\mu|
=(1/q)^{\nu_\gerX(P)}$. We will use the following lemma in the
next subsection.

\begin{lemma}{\label{lemma:pullback-frobenius}}
The morphism of sheaves $\pr_{\tilde{\gamma}_P}^*$ reduces to zero
modulo $\varpi/\mu$.

\end{lemma}

\begin{proof}  We will denote reduction modulo $\varpi/\mu$ of a
morphism of the form $\tilde{\gamma}_Q$ by $\bar{\gamma}_Q$. For
simplicity, denote the $\kappa$-algebra $\calO_E/(\varpi/\mu)$ by
$R$. Let
\[
s^\prime:X \ts \kappa \ts_{\kappa} R \ra Y\ts \kappa \ts_{\kappa}
R
\]
be the base extension from $\kappa$ to $R$ of the section $s$
given in assumption {\bf A2.1}. By Proposition 3.10 of
\cite{GorKas} we have $\bar{\gamma}_{\gers_\rig(P)}=s^\prime
\bar{\gamma}_P$. Hence the morphism
$\pr_{\tilde{\gamma}_P}^*=\tilde{\gamma}_{\gers_\rig(P)}^*\pr^*$
reduces modulo $\varpi/\mu$ to $\bar{\gamma}_P^*(s^\prime)^*
(\pr^*\ts_{\calO_{L_0}} R)$. But $(s^{\prime})^*
(\pr^*\ts_{\calO_{L_0}} R)$ is the base extension from $\kappa$ to $R$ of
\[
s^* (\pr^* \ts \kappa)
\]
which is zero by assumption {\bf A8}.
\end{proof}

\

\subsection{The main theorem.}  We now start preparing for the proof of the classicality result.
By the last assertion in assumption
{\bf{H3.1}},  we have a product diagram

\begin{eqnarray}\label{diagram:product-higherlevel}
\xymatrix{
(\phi_{1,\rig}\lambda_{1,\rig})^{-1}(\gerY_\rig(e/(e+1),1])
\ar[rr]^{\lambda_{1,\rig}} \ar[d]_{\eta_{1,\rig}} &&
(\phi_{1,\rig})^{-1}(\gerY_\rig(e/(e+1),1])  \ar[d]^{\phi_{1,\rig}}  \\
(\pi_{1,\rig})^{-1}(\gerY_\rig(e/(e+1),1]) \ar[rr]_{\pi_{1,\rig}}
&& \gerY_\rig(e/(e+1),1]}
\end{eqnarray}
By base extension the section $\gers_\rig^0$ defined in assumption
{\bf A6} yields a section to $\lambda_{1,\rig}$ defined on $(\phi_{1,\rig})^{-1}(\gerY_\rig(e/(e+1),1]) \subset \gerXm_\rig$ as follows
\[
\gert_\rig: (\phi_{1,\rig})^{-1}(\gerY_\rig(e/(e+1),1]) \
\rightarrow\ \gerXmo_\rig.
\]
 For any  subinterval $I$ of $(e/(e+1),1]$ we set
\[
\calV I:=\phi^{-1}_{1,\rig}(\gerY_\rig I) \subseteq (\phi_{1,\rig})^{-1}(\gerY_\rig(e/(e+1),1])
\]
Recall that for any interval $I$ as above we have set $I^\tau:=\{1-e(1-r):r\in I \}.$

\begin{lemma}\label{lemma:pre-tau}
For each interval $I \subseteq (e/(e+1),1]$, we have
$\gert_\rig(\calV I) \subseteq \lambda_{2,\rig}^{-1}(\calV
I^{\tau})$. In fact we have
\[
\gert_{\rig}  (\calV I) =\lambda_{2,\rig}^{-1}(\calV I^{\tau}) \cap \lambda_{1,\rig}^{-1}(\calV I).
\]
\end{lemma}

\begin{proof} This follows from the corresponding result over $\gerY_\rig$ (Corollary \ref{corollary:section-tame}) in the usual way. By definition of $\gert_{\rig}$ we can write
\begin{eqnarray*}
\gert_\rig(\phi^{-1}_{1,\rig}(\gerY_\rig I))=
\eta_{1,\rig}^{-1}(\gers^{0}_{\rig}(\gerY_{\rig}I) )&=&
\eta_{1,\rig}^{-1}\pi_{2,\rig}^{-1}(\gerY_{\rig} I^{\tau}) \cap \eta_{1,\rig}^{-1}\pi_{1,\rig}^{-1} (\gerY_{\rig} I)  \\
&=&\lambda_{2,\rig}^{-1}\phi_{1,\rig}^{-1}(\gerY_{\rig}I^{\tau}) \cap \lambda_{1,\rig}^{-1}\phi_{1,\rig}^{-1}(\gerY_{\rig}I) \\
&=&\lambda_{2,\rig}^{-1}(\calV I^{\tau}) \cap \lambda_{1,\rig}^{-1}(\calV I).
\end{eqnarray*}
\end{proof}

\

\begin{defn}{\label{definition:tau}} Let $\tau\colon \calV(e/(e+1),1]
\ra\gerXm_\rig$ be given by $\tau=\lambda_2\gert_\rig$. By Lemma
\ref{lemma:pre-tau}, for any interval $I\subseteq (e/(e+1),1]$, we
have a morphism
\[
\tau: \calV I  \ra \calV I^{\tau}.
\]
We have
\[
\tau^*\omg=\gert_\rig^*\lambda_{2,\rig}^*\omg=\gert_\rig^*\lambda_{2,\rig}^*\phi_{1,\rig}^*\omg=\gert_\rig^*\eta_{1,\rig}^*\pi_{2,\rig}^*\omg=\gert_\rig^*\eta_{1,\rig}^*(\pi_{1,\rig}^\prime)^*w_{\rig}^*\pi_\rig^*\omg.
\]
Using this and commutative diagrams given by assumption {\bf H3.1}, we see that  applying
$\gert_\rig^*\eta_{1,\rig}^*(\pi_{1,\rig}^\prime)^*$ to the
morphism $\pr^*$ given by Definition \ref{definition:omg} yields a morphism of
$\calO_{\calV I}$-modules
\[
\pr^*:=\gert_\rig^*\eta_{1,\rig}^*(\pi_{1,\rig}^\prime)^*\pr^*:\tau^*\omg\ra \omg
\]
Comparing with the definition of $\vartheta$ given right before Lemma \ref{lemma:redundant} it is clear that $(\pr^*)^{\ts k}$ is obtained by specializing $\vartheta$ via $\gert_\rig$.

Remembering $\gerF_\rig=\omg^{\ts k}$ on $\gerXm_\rig$, define $\tau:H^0(\calV I^{\tau}, \gerF_\rig) \ra H^0(\calV I,
\gerF_\rig)$ by
\[
f \mapsto f^\tau:=\varpi^{-k}(\pr^*)^{\ts k}\tau^*f.
\]
\end{defn}

\

\begin{rem} \label{remark:tau-sigma} In particular, taking $I=[1-r,1-r]$ for $r<1/(e+1)$, we find that if $\nu_\gerY(\phi_{1,\rig}(Q))=1-r$ then $\nu_\gerY(\phi_{1,\rig}(\tau(Q)))=1-er$.  If we assume further that $r<1/e(e+1)$ then using Proposition \ref{proposition:GK} (2) we find that if $R_0:=\pi_{\rig} (\phi_{1,\rig}(Q))$ and $R_1:=\pi_{\rig} (\phi_{1,\rig}(\tau(Q)))$ then $\nu_\gerX(R_1)=e\nu_\gerX(R_0)=e^2r$. We will use this later.
\end{rem}

 There is an alternative way to describe $\tau:H^0(\calV I^{\tau}, \gerF_\rig) \ra H^0(\calV I,
\gerF_\rig)$ using  the machinery of Definition \ref{definition:trace}.
We will use this description in the  proof of Proposition
\ref{proposition:extension}.

Let the data $(\calG,\calY^0,\calY,\calV,\calW,\alpha_1,\alpha_2,\ell)$ in Definition \ref{definition:trace} be given by
\[
(\gerF_\rig,\gerXmo_\rig,\gerXm_\rig,\calV I^\tau \cup \calV I^w, \calV I, \lambda_{1,\rig},\lambda_{2,\rig},\vartheta)
\]
To do so we need to check $\lambda_{1,\rig}^{-1}(\calV I)\subseteq
\lambda_{2,\rig}^{-1}(\calV  I^{\tau}) \cup \lambda_{2,\rig}^{-1}( \calV I^{w})$ which follows from Lemma \ref{lemma:pre-trace} by our usual trick of applying $\eta_{1,\rig}^{-1}$ and using  the commutative diagrams in assumption {\bf H3.1} for $i=1$. Also note that this is an admissible disjoint union by Lemma \ref{lemma:pre-trace}. We get  a map
\[
T_{\calV I^\tau \cup \calV I^w}^{\calV I}: H^0(\calV I^\tau, \gerF_\rig) \oplus H^0(\calV I^w,\gerF_\rig)  \ra H^0(\calV I,\gerF_\rig)
\]
\begin{lemma}{\label{lemma:alternate-tau}}
For any $f \in H^0(\calV I^{\tau}, \gerF_\rig)$ we have $f^\tau=\varpi^{-k}T_{\calV I^\tau \cup \calV I^w}^{\calV I}(f,0)$.
\end{lemma}

\begin{proof} By definition (and since $(f,0)$ is identically zero on $\calV I^w$), we have
\[
T_{\calV I^\tau \cup \calV I^w}^{\calV I}(f,0)=(\lambda_{1,\rig})_*(\vartheta \lambda_{2,\rig}^*f)_{|_\calD}
\]
where $\calD=\lambda_{2,\rig}^{-1}(\calV I^\tau) \cap \lambda_{1,\rig}^{-1} (\calV I)$ and $\lambda_{2,\rig}$ is considered only as a map from $\lambda_{2,\rig}^{-1}(\calV I^\tau)$ to $\calV I^\tau$.
By Lemma \ref{lemma:pre-tau} we have $\calD=\gert_\rig(\calV I)$ and hence $(\lambda_{1,\rig})_*$ can be rewritten as $\gert_\rig^*$. To finish the proof we note that $\gert_\rig^* \vartheta=(\pr^*)^{\ts k} \gert_\rig^*$ as was noted in Definition \ref{definition:tau}.

\end{proof}

\begin{prop}{\label{proposition:tau-norm}}
Let $I$ be a closed interval in $(e/(e+1),1]$. Let $h \in H^0(\calV I^{\tau}, \gerF_\rig)$. For
any $Q \in \calV I$ we have
\[
|h^\tau(Q)| \leq q^{k\nu_\gerX(\pi_\rig\phi_{1,\rig}(Q))}
|h(\tau(Q))|
\]
where the norms are as in Definition \ref{definition:norm}.
\end{prop}

\begin{proof}
In this proof the notation $\gamma_.,\tilde{\gamma}_.$ is as in the
paragraph before Definition \ref{definition:norm}.   Let us fix
$\gamma_Q:\Sp(E) \ra \calV I$ giving the point $Q$, and
$\gamma_{\tau(Q)}:=\tau \gamma_Q$ give the point $\tau(Q)$. Denote
by $\tilde{\gamma}_Q$, and $\tilde{\gamma}_{\tau(Q)}$, respectively,
their formal extensions. Let $\tilde{\gamma}_{\gert_\rig(Q)}$ denote
the formal extension of $\gamma_{\gert_\rig(Q)}:=\gert_\rig
\gamma_Q$.

Let $\pr_{\gamma_Q}^*:\gamma_{\tau(Q)}^*\omg=\gamma_Q^*\tau^*\omg\ra \gamma_Q^*\omg
$ denote the $E$-linear morphism obtained by
specializing $\pr^*$ (as in Definition \ref{definition:tau}) via
$\gamma_Q$. In other words
$\pr_{\gamma_Q}^*=\gamma_{\gert_\rig(Q)}^*\eta_{1,\rig}^*(\pi_{1,\rig}^\prime)^*\pr^*$,
where now $\pr^*$ is as in assumption {\bf A7}. This has a formal
model
\[
\pr_{\tilde{\gamma}_Q}^*:=\tilde{\gamma}_{\gert_\rig(Q)}^*\eta^*(\pi_1^\prime)^*\pr^*:\tilde{\gamma}_{\tau(Q)}^*\omg \ra \tilde{\gamma}_Q^*\omg
\]
These morphisms fit into the following commutative diagram.

\begin{eqnarray}
\xymatrix{  H^0(\calV I^\tau, \omg^{\ts k}) \ar[r]^{\tau^*}  \ar[dr]_{\gamma_{\tau(Q)}^*} &  H^0(\calV I, \tau^* \omg^{\ts k}) \ar[rr]^{(\pr^*)^{\ts k}} \ar[d]_{\gamma_Q^*} && H^0(\calV I,  \omg^{\ts k}) \ar[d]_{\gamma_Q^*} \\
                      &  H^0(\Sp(E), \gamma_{\tau(Q)}^* \omg^{\ts k}) \ar[rr]^{(\pr_{\gamma_Q}^*)^{\ts k}} && H^0(\Sp(E), \gamma_Q^* \omg^{\ts k}) \\
                      &  H^0(\Spf(\calO_E), \tilde{\gamma}_{\tau(Q)}^*  \omg^{\ts k}) \ar[u]\ar[rr]^{ (\pr_{\tilde{\gamma}_Q}^*) ^{\ts k}}  & & H^0(\Spf(\calO_E), \tilde{\gamma}_Q^*  \omg^{\ts k})\ar[u]}
                      \end{eqnarray}
Let $\mu \in \calO_E$ be such that $|\mu|
=(1/q)^{\nu_\gerX(\pi_\rig\phi_{1,\rig}(Q))}$. To prove the
statement, is enough to show that $|h(\tau(Q))|\leq 1$ implies
$|h^\tau(Q)| \leq  |1/\mu|^k$.  We have
\[
|h^\tau(Q)|=|\gamma_Q^*h^\tau|_Q=|\varpi^{-k} \gamma_Q^*
(\pr^*)^{\ts k} \tau^*h|_Q=|\varpi|^{-k} |(\pr_{\gamma_Q}^*)^{\ts
k}\ \gamma_{\tau(Q)}^*  h|_Q
\]
Unraveling the definitions of the norms shows that it suffices to
prove
\[
\pr_{\gamma_Q}^* {\large
(}H^0(\Spf(\calO_E),\tilde{\gamma}_{\tau(Q)}^*\omg) {\large )}
\subseteq (\varpi/\mu) H^0(\Spf(\calO_E), \tilde{\gamma}_Q^* \omg
).
\]
Since $\pr_{\tilde{\gamma}_Q}^*$ is the restriction of
$\pr_{\gamma_Q}^*$ to $\tilde{\gamma}_{\tau(Q)}^*\omg \subset
\gamma_{\tau(Q)}^* \omg$, it is enough to show that the
reduction of $\pr_{\tilde{\gamma}_Q}^*$ modulo $\varpi/\mu$ is the
zero morphism. Let $R=\phi_{1,\rig}(Q)$, and $P=\pi_\rig(R)$.
Then, we have
\[
\pr_{\tilde{\gamma}_Q}^*=\tilde{\gamma}_{\gert_\rig(Q)}^*\eta^*(\pi_1^\prime)^*\pr^*
=\phi_1^*\tilde{\gamma}_{\gers^0_\rig(R)}^*(\pi_1^\prime)^*\pr^*=
\phi_1^*\pi^*\tilde{\gamma}_{\gers_\rig(P)}^*\pr^*=
\phi_1^*\pi^*\pr_{\tilde{\gamma}_{\sigma(P)}}^*,
\]
where $\tilde{\gamma}_{\gers^0_\rig(R)}$ (respectively,
$\tilde{\gamma}_{\gers_\rig(P)}$) denotes the formal extension of
$\gamma_{\gers^0_\rig(R)}:=\gers^0_\rig\gamma_R$ (respectively,
$\gamma_{\gers_\rig(P)}:=\gers_\rig \gamma_P$). Note that the second
(respectively, third) equality comes from the fact that $\gert_\rig$
(respectively, $\gers^0_\rig$) is obtained from $\gers^0_\rig$
(respectively, $\gers_\rig$) by base extension. Now the result
follows since by Lemma \ref{lemma:pullback-frobenius} we know that
$\pr_{\tilde{\gamma}_{\sigma(P)}}^*$ reduces to the zero morphism
modulo $\varpi/\mu$.

\end{proof}

\begin{cor}\label{corollary:ordbnd}
Let~$h \in H^0(\calV[1,1],\gerF_\rig)$. For all~$n \geq 0$ we have
 $|h^{\tau^n}|_{_{\calV[1,1] }} \leq |h|_{_{\calV[1,1] }} < \infty$.
 \end{cor}

\begin{proof}
This follows from Proposition \ref{proposition:tau-norm} with $I=[1,1]$, and Lemma
\ref{lemma:Banach}. We only remind the reader that $\pi_\rig(\gerY[1,1])=\gerX[0,0]$.
\end{proof}

For a generalized eigenform $f$ of $\U$ with eigenvalue $a$ we
define the slope of $f$ to be $\val(a)$. Recall our fixed choice of
(the ``normalization factor'')  $c \in L_0$ in Definition
\ref{definition:trace}. It was used to define the various $\U$ operators (See Remark \ref{remark:after-trace}). Also recall that we have fixed a weight $k$ and a choice of $\pr^*$ satisfying Assumption {\bf A8} for $\calF$.

\begin{thm}{\label{theorem:classicality}} Let either $0<r<e^{2-m}/(e+1)$  and $f\in S_r(\gerXm_\rig,\gerF_\rig;L)$, or $0<r<e/(e+1)$ and $f\in S_r(\gerY_\rig,\gerF_\rig;L)$.
\begin{enumerate}

\item Let $A(x) \in L[X]$ be a polynomial such that all roots of
$A$ in the algebraic closure of $L$ have valuation less than
$k+\val(c)$. If $A(\U)f$  is classical, then so is $f$.

\

\item (Special case of the first part). If $f$ is a generalized
$\U$-eigensection  which has a weight $k$ and slope less than $k+\val(c)$,
then $f$ is classical.

\end{enumerate}

\end{thm}

\begin{proof} Without loss of generality we can assume $L=L_0$. We remark that it is enough to prove the result over $\gerX(m)_\rig$ as the proof for the other case follows exactly in the same way, or alternatively by observing that $\gerY_\rig$  satisfies the axioms required to be an instance of $\gerX(m)_\rig$ for $m=1$ with obvious maps.
It is also enough to deal with the case when $A$ has degree one. Assume we have done so. Then the general case can be proved  by an induction as follows. Assume $A(\U)f=F$ is classical. Passing to a finite extension of $L$ we can assume~$A(x)=(x-a_1)(x-a_2)...(x-a_l)$, such that~$\val(a_j)<k+\val(c)$ for
all~$j$.  For $1 \leq j \leq l-1$ define~$f_j=(\U-a_{j+1})(\U-a_{j+2})...(\U-a_l)f$. Then we
have~$(\U-a_1)f_1=F$ and hence by assumption~$f_1$ is classical.
Similarly we have~$(\U-a_2)f_2=f_1$ and since~$f_1$ is classical,
we deduce that~$f_2$ is classical. Continuing this way we see
that $(\U-a_l)f$ is classical  and hence, by assumption, $f$ is classical.

Assume now  that $\U f-af$ can be extended to $F \in
H^0(\gerXm_\rig,\gerF_\rig)$, and that $\val(a)<k+\val(c)$. For any interval of the form $I=[x,1)$, we let $\bar{I}=[x,1]$. Let us
fix a rational $0<r<1/(e+1)$. Let $I_n=[1-re^{-n},1)$ for $n \geq 0$.
Then $I_0 \supset I_1 \supset \cdots $ and we have $I_{n+1}^\tau=I_n$.

Let us denote $\phi_{1,\rig}^{-1}(\gerY_\rig[0,1))$ by $\calU$ for simplicity.  By Proposition \ref{proposition:analytic-higher} $f$ can be extended to $\calU$, and we have $\U
f=af+F_{|_{\calU}}$, where $\U$ denotes $\U_\calU^\calU=cT_\calU^\calU$ in the
notation of Definition \ref{definition:trace}. Recall the definition of $\calV I$ from \S 4.2. We will denote the
restriction of $f,F$ to $\calV I_n$ by the same letters. Let $b:=c\varpi^k/a$. We have $\val(b)>0$.
\begin{prop} {\label{proposition:extension}}The section $f-bf^\tau  \in H^0( \calV I_1,\gerF_\rig )$ extends to a section $F_1$ in
$H^0(\calV {\bar I_1},\gerF_\rig)$.
\end{prop}

\begin{proof} Recall from the discussion leading to Lemma \ref{lemma:alternate-tau} the operator $T_{\calV I_1^\tau \cup \calV I_1^w}^{\calV I_1}$. We can write
\[
(\U f)_{|_{\calV I_1}}=c(T_\calU^\calU f)_{|_{\calV I_1}}=
cT_{\calV I_1^\tau \cup \calV I_1^w}^{\calV I_1} (f,f)
\]
where $(f,f) \in H^0(\calV I_1^\tau,\gerF_\rig) \oplus H^0(\calV I_1^w,\gerF_\rig)$ and for the second equality we have used Lemma \ref{lemma:properties-of-trace} (1). But we can write
\[
T_{\calV I_1^\tau \cup \calV I_1^w}^{\calV I_1} (f,f)=T_{\calV I_1^\tau \cup \calV I_1^w}^{\calV I_1} (f,0)+T_{\calV I_1^\tau \cup \calV I_1^w}^{\calV I_1} (0,f)=\varpi^kf^\tau+T_{\calV I_1^\tau \cup \calV I_1^w}^{\calV I_1} (0,f)
\]
using Lemma \ref{lemma:alternate-tau}. Therefore, we have the following equation of sections of $\gerF_\rig$ on $\calV I_1$
\[
f-bf^\tau=(c/a)T_{\calV I_1^\tau \cup \calV I_1^w}^{\calV I_1} (0,f)-F/a,
\]
and to prove the result it is enough to show that $T_{\calV I_1^\tau \cup \calV I_1^w}^{\calV I_1} (0,f)$ can be extended to $\calV \bar{I}_1$.
But such an extension is provided by $T_{\calV {\bar I}_1^\tau \cup \calV {\bar I}_1^w}^{\calV {\bar I}_1} (0,f) $ which is well defined by the discussion before Lemma \ref{lemma:alternate-tau}
and the fact that $f$ is indeed defined on $ \calV {\bar I}_1^w \subseteq  \phi_{1,\rig}^{-1}(\gerY_\rig[0,1/(e+1)])$.
\end{proof}

\begin{lemma}{\label{lemma:f-bounded}}
We have $|f|_{\calV I_0}< \infty$.
\end{lemma}

\begin{proof}
By Proposition \ref{proposition:extension} $f-bf^\tau$ extends from $\calV I_1$ to $\calV {\bar{I_1}}$  which is an affinoid by Lemma \ref{lemma:affinoid}.  Hence by Lemma \ref{lemma:Banach} $f-bf^\tau$ will have finite norm on $\calV I_1$. Let $M^{\prime\prime}$ be a common upper  bound for this norm and $|f|_{_{\calV[1-r,1-re^{-1}]}}$.   We prove by induction that for all $n \geq 0$, $f$ is bounded on $\calV_n\!:=\!\calV[1-re^{-n},1-re^{-n-1}]$ by $M_n:=M^{\prime\prime} q^{(kr +\cdots +kre^{-n+1})}$ (we let $M_0=M^{\prime\prime}$). The result will then follow as $\calV I_0$ is the union of $\calV_n$'s and $M_n \leq M^\prime:=M^{\prime\prime} q^{kre/(e-1)}$.

The claim is true for $n=0$. Let~$Q \in \calV_{n+1}$. Then~$\tau(Q) \in \calV_n$. Let $P:=\pi_\rig(\phi_{1,\rig}(Q))$. Since $\nu_\gerY(\phi_{1,\rig} Q) >e/(e+1)$ part (2) of Proposition \ref{proposition:GK} tells us  $\nu_\gerX(P)=e(1-\nu_\gerY(\phi_{1,\rig}(Q)))\leq re^{-n}$ . By  Proposition \ref{proposition:tau-norm} and  the
induction hypothesis we have
\[
|f^\tau(Q)| \leq q^{k\nu_\gerX(P)}|f(\tau(Q))|
\leq  q^{kre^{-n}}M_n = M_{n+1 }
\]
which implies that $|f^{\tau}|_{_{ \calV_{n+1}}} \leq M_{n+1}$. So we can write
\[
|f|_{_{\calV_{n+1}}} \leq \max\{ |f-bf^\tau|_{_{\calV_{n+1}}} , |b{f^\tau}|_{_{\calV_{n+1}}}  \}
\leq \max\{M^{\prime\prime}, |{f^\tau}|_{_{\calV_{n+1}}} \}
\leq M_{n+1}.
\]

\end{proof}

\begin{lemma}{\label{lemma:uniform-bdd}}
There is an $M >0$ such that for all $n \geq 0$ we have $|f^{\tau^{n}}|_{_{\calV I_{n}}} \leq M$.
\end{lemma}

\begin{proof}
Let~$Q \in \calV I_{n}$ with $n \geq 1$. By Proposition \ref{proposition:tau-norm}, we have $|f^{\tau^n}(Q)| \leq q^{k\nu_{\gerX}(P_0)}|f^{\tau^{n-1}}(\tau(Q))|$ where  $P_0$ denotes $\pi_\rig(\phi_{1,\rig}(Q))$. Note that by part (2) of Proposition \ref{proposition:GK} $\nu_\gerX(P_0)=e(1-\nu_\gerY(\phi_{1,\rig}(Q))) \leq re^{-n+1}$. Inductively, we find
\[
|f^{\tau^n}(Q)| \leq |f(\tau^n(Q))|
  \prod_{j=0}^{n-1} q^{k\nu_{\gerX}(P_j)}
\]
where $P_j=\pi_\rig(\phi_{1,\rig}(\tau^j(Q)))$. By Remark \ref{remark:tau-sigma} we have $\nu_\gerX(P_j)=e^j\nu_\gerX(P_0) \leq re^{j-n+1}$. Also since $\tau^n(Q) \in \calV I_0$ Lemma \ref{lemma:f-bounded} gives us $|f(\tau^n(Q))| \leq M^{\prime}$. Together, these imply

\begin{eqnarray}
|f^{\tau^n}(Q)| \leq M^{\prime}  q^{kr\Sigma_{j=0}^{n-1} e^{j-n+1}}\leq   M^{\prime}q^{kre/(e-1)}=:M.\nonumber
\end{eqnarray}

\end{proof}

Next we recall a gluing lemma which was proved in \cite{Kassaei3}.

\begin{lemma}{\label{lemma:gluing}}
Let $\gerZ$ be a quasi-compact
admissible formal scheme over $\calO_0$ and $\gerN$ an invertible
sheaf on it.  Let~$\calX \subset \gerZ_\rig$ be a
smooth affinoid subdomain. Assume that~$\calX$ is a disjoint union
of two admissible opens~$\calX=\calY \cup \calZ$, where $\calZ$ is an
affinoid. Assume we are given affinoid subdomains of $\calX$ denoted
by~$\calZ_n$ for~$n \geq 1$ with
\begin{eqnarray}{\label{equation:corrected}}
\calZ \subset \cdots  \calZ_3 \subset \calZ_2 \subset  \calZ_1
\end{eqnarray}
and such that~$\{\calY,\calZ_n\}$ is an admissible cover of~$\calX$ for each
$n$. Assume that we are given two sections
\[
{\bf f}\in H^0(\calY,\gerN_\rig) \ \ \ \ \ \ \ \ \ {\bf g}\in H^0(\calZ,\gerN_\rig)
\]
and for each~$n \geq 1$, a section~$F_n \in H^0(\calZ_n,\gerN_\rig)$
such that, as~$n \ra \infty$, we have
\[
|F_n-{\bf f}|_{_{\calY \cap \calZ_n}} \ra 0\ \ \ {\rm and}\ \ \
|F_n-{\bf g}|_{_\calZ} \ra 0.
\]
Then~${\bf f}$ and~${\bf g}$ glue together to give a global section
of~$\gerN_\rig$ on~$\calX$. In other words, there is a section
of~$\gerN_\rig$ on~$\calX$, which restricts to~${\bf f}$ on~$\calY$, and
restricts to~${\bf g}$ on~$\calZ$.
\end{lemma}

\begin{rem}{\label{remark:shu}}
In \cite{Kassaei3} all but the first of the inclusions in Equation (\ref{equation:corrected}) are written in the reverse order. We are thankful to Shu Sasaki for pointing out this typo. The proof in \cite{Kassaei3} is written with the correct inclusions in mind!

\end{rem}
\end{proof}

We want to use this lemma to glue $f$ on $\calU=\phi_{1,\rig}^{-1}\gerY[0,1)$ with a section $g$  that we will construct below on $\phi_{1,\rig}^{-1}\gerY[1,1]=\calV[1,1]$ and produce a classical section of $\gerF_\rig$ on $\gerXm_\rig$.

Consider  the section $F_1$ constructed in Proposition \ref{proposition:extension}.
Define  $F_n:=\sum_{i=0}^{n-1} b^i F_1^{\tau^i} $ in $H^0(\calV {\bar I}_n,\gerF_\rig)$.  Since $F_1$ restricts to $f-bf^\tau$ over $\calV I_1$, we easily see that ${F_n}_{|_{\calV I_n}}=f-b^nf^{\tau^n}$. Since $|b|<1$ and by Corollary \ref{corollary:ordbnd} and Lemma \ref{lemma:Banach} we can define $g:=\sum_{n=0}^\infty b^i F_1^{\tau^i} $ in $H^0(\calV [1,1],\gerF_\rig)$.

We want to apply the gluing lemma
with~$\gerZ=\gerXm$,~$\gerN=\gerF_\rig$,~$\calX=\calV {\bar I}_0$ which is a smooth affinoid  by  assumption {\bf H1} and Lemma \ref{lemma:affinoid}, $\calY=\calV I_0$, $\calZ=\calV[1,1]$,~$\calZ_n=\calV \bar{I}_n$,~${\bf f}=f$,~ and ${\bf g}=g$. We
have~$F_n-f=-b^nf^{\tau^n}$ on $\calV I_n=\calY \cap \calZ_n$ and
$F_n-g=-b^ng^{\tau^n}$ on $\calV[1,1]$, and  therefore  in view of Lemma \ref{lemma:uniform-bdd} and Corollary \ref{corollary:ordbnd}, and since
$|b|<1$, we can apply Lemma \ref{lemma:gluing} to obtain a  section of $\gerF_\rig$ on $\calV {\bar I}_0=\calV[1-r,1]$ denoted $f^\prime$ for the moment.  By construction $f^\prime$ and $f \in H^0(\calU,\gerF_\rig)$ restrict to the same section on $\calV I_0$. Since $\{
\calV {\bar I}_0,\calU \}$ forms an admissible covering of $\gerXm_\rig$ we obtain a global section of $\gerF_\rig$ on $\gerXm_\rig$ which extends $f$, and hence $f$ is classical.

\section{applications}
As we were setting up the notation and progressing in the first part of this paper, we explained how the case of modular curves (where one has to take the normalization factor $c$ from Definition \ref{definition:trace} to be $1/p$) is an example covered by our results. In this section we show how the case of various Shimura curves are also covered.

 \subsection{Unitary Shimura curves}

In \cite{Kassaei2}  we developed a theory of overconvergent modular forms over certain unitary Shimura curves, and stated that we expected the analogue of Coleman's {\it classicality} result (cf. \cite{Coleman1,Coleman2}) that ``overconvergent modular forms of small slope are classical''  to be true over these Shimura curves. We now show how this follows as a special case of Theorem \ref{theorem:classicality}. In fact we do  more: in \cite{Kassaei2} we only studied overconvergent modular forms of level $\Gamma_0(\mathcal{P})$ (which in the notation of this paper corresponds to $\gerY_\rig$), whereas with results that we have proven for $\gerXm_\rig$ (the ``higher-level" cases) we can now extend the constructions of \cite{Kassaei2} to  the case where level structures contain arbitrary powers of $\mathcal{P}$ (we will make this precise below),  and provide a classicality result for these overconvergent modular forms of higher levels as well.  For simplicity of referencing we will stay faithful to the notation of \cite{Kassaei2} to a large extent, even though at times it may not be the most economical one.

Let $F$ be a totally real field of degree $d>1$. Let $\calP_1=\calP,\calP_2,\cdots,\calP_r$ be the primes of $F$ over $p$. Let $F_{\calP_i}$ be the completion of $F$ at $\calP_i$ with ring of integers $\calO_{\calP_i}$. Let $B$ be a quaternion algebra over $F$ which splits at $\calP$ and also at exactly one infinite place of $F$. Choose $\lambda<0$ a rational number such that $\QQ(\lambda)$ splits at $p$ and define $E=F(\lambda)$. It follows that the primes of $E$ lying above $p$ come in pairs, each pair lying over one of the $\calP_i$'s, and one gets an isomorphism
\begin{eqnarray}\label{equation:decompose}
\calO_E\ts \ZZ_p   \eqra (\calO_{{\mathcal P}_1}\oplus \cdots \oplus
\calO_{{\mathcal P}_r})\oplus (\calO_{{\mathcal P}_1}\oplus \cdots \oplus \calO_{{\mathcal
P}_r})
\end{eqnarray}
Define $D:=B \ts_F E$ and let $V$ denote $D$ as a $\QQ$-vector space
with a left action of $D$. Let  $\calO_D=:V_\ZZ \subset V=D$ be a
maximal order of $D$.  The  isomorphism \ref{equation:decompose}
induces the following decompositions of $D\ts \QQ_p$ and $\calO_D
\ts \ZZ_p$
\[
\begin{array}{ccccccccccccc} \calO_D \ts
\ZZ_p&=&\calO_{D^1_1}&\oplus&\cdots&\oplus&\calO_{D^1_r}&\oplus&\calO_{D^2_1}&\oplus&\cdots&\oplus&\calO_{D^2_r}\\
\bigcap&   &\bigcap& & &&\bigcap&  & \bigcap& & & &\bigcap\\
                  D \ts \QQ_p&=&D^1_1&\oplus&\cdots&\oplus &D^1_r &\oplus& D^2_1&\oplus&
\cdots&\oplus& D^2_r \end{array}
\]
where each $D^k_j$ is an
$F_{\calP_j}$-algebra isomorphic to $B\ts_F F_{\calP_j}$. In
particular, $D^1_1$ and $D^2_1$ are isomorphic to ${\rm M}_2(F_{\calP})$.

One can choose $\calO_D$, an involution of second type $l \mapsto
l^*$ on $D$, and a $\QQ$-valued alternating non-degenerate  form
$\Psi$ on $V$   satisfying $\Psi(lv,w)=\Psi(v,l^*w)$ for $v,w \in V$
and $l \in D$ such that
  \begin{itemize}
\item[i)] $\calO_D$ is stable under the involution $l \mapsto l^*$; in fact the involution switches $D^1_j$ and $D^2_j$.
\item[ii)] Each $\calO_{D^k_j}$ is a maximal order in $D^k_j$ and
$\calO_{D^2_1}\subset D^2_1={\rm M}_2(\calF_\calP)$ is identified with ${\rm M}_2(\calO_\calP)$,
\item[iii)] $\Psi$ takes integer values on $V_\ZZ$,
\item[iv)] $\Psi$ induces a perfect pairing $\Psi_p$ on
$V_{\ZZ_p}=V_\ZZ \ts \ZZ_p.$
\end{itemize}
Each $\calO_D\ts\ZZ_p$-module (or any element of an abelian category with an action by $\calO_D\ts\ZZ_p$) $\Lambda$, hence, admits a decomposition
\begin{eqnarray}{\label{equation:decomposition}}
\Lambda=\Lambda^1_1\oplus...\oplus\Lambda^1_r\oplus\Lambda^2_1\oplus...\oplus\Lambda^2_r
\end{eqnarray}
such
that each $\Lambda^k_j$ is an $\calO_{D^k_j}$-module. The
${\rm M}_2(\calO_\calP)$-module $\Lambda^2_1$ can be further decomposed as the
direct sum of two $\calO_\calP$-modules $\Lambda^{2,1}_1$ and
$\Lambda^{2,2}_1$ by choosing idempotents  in ${\rm M}_2(\calO_\calP)$.

Let $G^\prime$ be the algebraic group which for any $\QQ$-algebra $R$ has  $R$-points given by the group of symplectic similitudes of $(V\ts_\QQ R,\Psi \ts_\QQ R)$.
The finite adelic points of $G^\prime$ can be described as
\[
G^\prime(\AA^\infty)=\QQ_p^\times \times \GL_2(\calF_\calP) \times (B\ts_F F_{\mathcal{P}_2})^\times
\times \cdots \times (B\ts_F F_{\mathcal{P}_r})^\times \times
G^\prime(\AA^{\infty,p}).
\]
We consider open compact subgroups $K^\prime$ of $G^\prime(\AA^\infty)$ of the form
\[
K^\prime=\ZZ_p^\times \times K^\prime_\calP \times H^\prime
\]
where $K^\prime_\calP$ is a subgroup of $\GL_2(F_\calP)$, $H^{\prime}$ is a subgroup of $ (B\ts_F F_{\mathcal{P}_2})^\times
\times \cdots \times (B\ts_F F_{\mathcal{P}_r})^\times \times
G^\prime(\AA^{\infty,p})
$, and such that $K^\prime$ is small enough so that it leaves stable the lattice $V_{\hat{\ZZ}}:=V_\ZZ \ts \hat{\ZZ} \subset V\ts \AA^\infty$.

The unitary Shimura curve $M^\prime_{K^\prime}$ defined over $F_\calP$ represents the functor from the category of $F_\calP$-schemes to the category of sets where any $S=\Spec(R)$ (where $R$ is an
$F_\calP$-algebra) is mapped to the set of isomorphism classes of all quadruples $(A,i,\theta,\bar{\alpha})$ such that

\begin{itemize}
\item[i)] $A$ is an abelian scheme of relative dimension $4d$ over $R$ with
an action of $\calO_D$ via $i: \calO_D \inc \End_R(A)$, which satisfies

\begin{itemize}
\item[1)] the projective $R$-module ${\Lie}_1^{2,1}(A)$ has rank $1$
and $O_\calP$ acts on it via $O_\calP \hookrightarrow R$,
\item[2)] for $j \geq 2$, we have ${\Lie}^2_j(A)=0$,
\end{itemize}

\item[ii)]$\theta$ is a polarization of $A$ (of degree prime to $p$) such that the
corresponding Rosati involution sends $i(l)$ to $i(l^*)$,

\item[iii)]$\bar{\alpha}$ is a $K^\prime$ level structure: it is a class modulo
$K^\prime$ of symplectic $\calO_D$-linear isomorphisms
$\alpha:{\hat{T}}(A) \eqra V_{\hat{\ZZ}}$.
\end{itemize}
Here ${\hat{T}}(A)=\prod_{l} T_l(A)$ denotes $\varprojlim_n A[n]$ as a sheaf
over $\Spec(R)$ in the \'etale topology and the symplectic form
on ${\hat{T}}(A)$  is the $\theta$-Weil pairing.  Also note that $\Lie(A)$ has an action of $\calO_D \ts \ZZ_p$ and $\Lie^{2,1}_1(A)$ and $\Lie^2_j(A)$ are defined as in (\ref{equation:decomposition}).

For any such abelian scheme $A$ we can consider various objects with an $(\calO_D \ts \ZZ_p)$-action and decompose them as in (\ref{equation:decomposition}). For example any $\calO_D$-invariant subgroup scheme $H$ of $A$ which is killed by a power of $q$ has an action of $\calO_D \ts \ZZ_p$. In particular $A[q^m]$ can be decomposed and $A[q^m]^{2,j}_1$ is defined for $j=1,2$ and has an action of $\calO_\calP$. We define $A[\varpi^m]^{2,j}_1$ to be the $\varpi^m$-torsion in  $A[p^m]^{2,j}_1$. It is an $\calO_\calP$-module scheme of rank $q^{2m}$. Note that $\theta$, being prime to $p$, induces an isomorphism when restricted to the $q^m$- torsion subgroup, and since the involution switches $D^1_j$ and $D^2_j$ we find that $\theta:A[q^m]^1_j \ra (A[q^m]^2_j)^\vee$ is an isomorphism. We refer the reader to \S 4.4 of \cite{Kassaei2} for the definition of a type 1 and type 2 subgroup scheme of $A$. The definition can be extended in an evident way to subgroups that are killed by  a power of $q$ (rather than just $q$). Given an $\calO_\calP$-invariant subgroup scheme $C$ of $A[\varpi^m]^{2,1}_1$, using the above duality, condition ii) in the definition of type 1,2 subgroups, and a conjugation between idempotents in $\calO_{D^2_1} \eqra M_2(\calO_\calP)$, the subgroup $C$ can be uniquely extended to subgroups of both type 1 and type 2 of $A$ which we denote respectively  by $t_1(C)$ and $t_2(C)$ both of rank $q^{m\dim(A)}$.

Let $\epsilon_A: A \ra \Spec(R)$ be the structure map. Then $\epsilon_{A,*}(\Omega_{A/R})$ has an action of $\calO_D \ts \ZZ_p$ and we define $\omg_{A/R}$ to be $\epsilon_{A,*}(\Omega_{A/R})^{2,1}_1$. The above conditions on $\Lie(A)$ show that $\omg_{A/R}$ is a line bundle on $\Spec(R)$. This construction can be done universally and will give us a line bundle $\omg$ on $M^\prime_{K^\prime}$.

For specific choices of $K^\prime_\calP \subset \GL_2(F_\calP)$ we will re-interpret a $K^\prime$ level structure. Let ${\hat T}^p(A):=\prod_{l\not =p}T_l(A)$ and
 denote $(T_p(A))^2_2\oplus \cdots \oplus
(T_p(A))^2_r$  by  $T^{\calP}_p$. Similarly let ${\hat W}^p:=V_\ZZ \ts {\hat \ZZ}^p$ and
denote by $W^{\calP}_p$   the direct sum $(V_{\ZZ_p})^2_2 \oplus
...\oplus(V_{\ZZ_p})^2_r$.

 If $K^\prime_\calP=K^\prime_0(\calP)$, i.e., the group of all  matrices in $\GL_2(\calO_\calP)$ whose left lower corner entry is congruent to $0$ modulo $\calP$, then a $K^\prime$ level structure can be thought of as a choice of $(C,\bar{\alpha}^{\calP})$ where
\begin{itemize}
\item[1)]  $C$ is a finite flat $\calO_\calP$-submodule scheme of rank $q$ of  $(A[\varpi])^{2,1}_1$;
\item[2)] $\bar{\alpha}^\calP$ is a class  of isomorphisms
$\alpha^{\calP}=\alpha^{\calP}_p \oplus \alpha^p: T^{\calP}_p(A) \oplus{\hat
T}^p(A) \eqra W^{\calP}_p \oplus {\hat W}^p$ modulo $H^\prime$, with $\alpha^{\calP}_p$
linear and $\alpha^p$ symplectic.
\end{itemize}
If $K^{\prime}_{\calP}=K^{\prime}_{1}(\calP^{m})$ consisting of matrices in $\GL_2(\calO_\calP)$ such that the upper and lower left corner entries are congruent to, respectively, $1$ and $0$ modulo $\varpi^{m}$, then a $K^\prime$ level structure can be
written as a choice of $(Q,\bar{\alpha}^{\calP})$ where
\begin{itemize}
\item[1)]  $Q$ is a point of exact $\calO_\calP$-order $\calP^m$ in $(A[\varpi^m])^{2,1}_1$,


\item[2)] $\bar{\alpha}^\calP$ is as above.
\end{itemize}
If $K^{\prime}_{\calP}=K^{\prime,0}(\calP)$ consisting of matrices in $K^\prime_0(\calP)$ whose upper right corner entry is divisible by $\varpi$, then  a $K^\prime$ level structure is a choice of
$(C,D,\bar{\alpha}^{\calP})$ where
\begin{itemize}
\item[1)] $C,\bar{\alpha}^\calP$ is as above;
\item[2)]  $D$ is a finite flat $\calO_\calP$-submodule scheme of rank $q$ of  $(A[\varpi])^{2,1}_1$ which intersects $C$ trivially.
\end{itemize}
Finally if $K^{\prime}_{\calP}=K^{\prime,0}_1(\calP^{m})$ consisting of matrices in  $K^{\prime}_{1}(\calP^{m})$ then a $K^\prime$ level structure is a choice of $(Q,D,\bar{\alpha}^{\calP})$ where
\begin{itemize}
\item[1)] $Q,\bar{\alpha}^\calP$ is as above;
\item[2)] $D$ is a finite flat $\calO_\calP$-submodule scheme of rank $q$ of  $(A[\varpi^m])^{2,1}_1$ which intersects the $\calO_\calP$-submodule scheme generated by $Q$ trivially.
\end{itemize}

When $K^{\prime}_{\calP}$ is, respectively, $K^{\prime}_{0}(\calP)$, $K^{\prime}_1(\calP^{m})$, $K^{\prime,0}(\calP)$, and $K^{\prime,0}_1(\calP^{m})$ we denote $M^{\prime}_{K^{\prime}}$ by, respectively,
 $M^{\prime}_{H^{\prime},0}(\calP)$, $M^{\prime}_{H^{\prime},1}(\calP^{m})$, $N^{\prime}_{H^{\prime},0}(\calP)$, $N^{\prime}_{H^{\prime},1}(\calP^{m})$. We will only consider integral models for the first two cases by considering the same moduli problem defined over $\calO_\calP$-algebras  where now $Q$ is a point of exact $\calO_\calP$-order $\calP^m$ in the  in the sense of Drinfeld: that is a map of $\calO_\calP$-modules $\phi:\calO_\calP/\calP^m \ra Hom(\Spec(R), (A[\varpi^m])^{2,1}_1)$ such that $\sum_{a \in \calO/\calP^m} [a]$ is a finite flat $\calO_\calP$-submodule scheme of rank $q^m$ of $(A[\varpi^m])^{2,1}_1$. Here $[a]$ is the closed subscheme  of $(A[\varpi^m])^{2,1}_1$ corresponding to the $R$-point $\phi(a)$ and by the sum of two closed subschemes we mean the closed subscheme given by the product of their ideals. We set $Q=\phi(1)$. Let us denote thes integral models by ${\bf M}^{\prime}_{H^{\prime},0}(\calP)$, ${\bf M}^{\prime}_{H^{\prime},1}(\calP^{m})$. We explain how these Shimura curves and maps between them provide examples for the set-up of the paper.

We set $L_{0}=F_{\calP}$, and $\calO_{0}=\calO_{\calP}$. First we
discuss the tame situation. Set $X={\bf
M}^{\prime}_{H^{\prime}}:={\bf
M}^{\prime}_{H^{\prime},1}(\calP^{0})$ a curve over $\calO_{\calP}$.
Set $Y={\bf M}^{\prime}_{H^{\prime},0}(\calP)$. Then $\gerX_{\rig}$
and $\gerY_{\rig}$ are, respectively, the $p$-adic analytifications
of $M^{\prime}_{H^{\prime}}:=M^{\prime}_{H^{\prime},1}(\calP^{0})$
and $M^{\prime}_{H^{\prime},0}(\calP)$. There is a morphism
$\pi:Y\ra X$  defined by forgetting $C$. The section $s: X \ts
\kappa \ra Y \ts \kappa$ is given by
$(A,i,\theta,\bar{\alpha}^\calP) \mapsto
(A,i,\theta,\bar{\alpha}^\calP,\Ker({\rm Frob}_q)^{2,1}_1)$. The
automorphism $w:Y \ra Y$ is defined by dividing a test object
$(A,i,\theta,C,\bar{\alpha}^\calP)$ by $t_2(C)$.  See \S 4.4 of
\cite{Kassaei2} for a precise definition of the quotient of
$(A,i,\theta,\bar{\alpha}^\calP)$ by $t_2(C)$. To that account we
only need to add the construction of the subgroup of order $q$, and
that will be given by $A[\varpi]^{2,1}_1/C$. The automorphism
$\delta:Y \ra Y$ multiplies $\alpha^\calP$ by $q^{-1}$. We can see
that the assumptions {\bf A1}-{\bf A3} hold in this case (except
that the slightly more relaxed version of {\bf A2.2} stated in
Remark \ref{finite extension} should be considered here since we are
taking $L_0=F_{\mathcal P}$) and $e=q$, either by using Carayol's
results in \cite{Carayol}, or by using the theory of local models
(see  \S 5 of \cite{GorKas} for a brief discussion). The results of
\cite{GorKas} apply to these curves. In particular there is a
measure of singularity $\nu_{\gerX}$ on $\gerX_{\rig}$ which we call
the measure of supersingularity in this case. And if
$\nu_{\gerX}(A,i,\theta,\bar{\alpha}^{\calP} )<q/(q+1)$ the image of
the section $\gers_{\rig}$ marks a unique subgroup scheme of order
$q$ in $A[\varpi]^{{2,1}}_{1}$. We call this the canonical subgroup
of $A[\varpi]^{{2,1}}_{1}$.
 Next we define $\gerY^{0}_{\rig}$ to be the analytification of $N^{\prime}_{H^{\prime},0}(\calP)$, which is a rigid analytic curve over $\calF_{\calP}$. There are two finite flat morphisms $\pi_{1,\rig}$ and $\pi^{\prime}_{1,\rig}$ from $\gerY^{0}_{\rig}$ to $\gerY_{\rig}$ which are, respectively, the analytifications of the morphisms forgetting $D$ and $C$. We need to verify {\bf A5} and {\bf A6}. The assumption {\bf A5} holds trivially, as $C$ and $D$ are different subgroups by assumption. The assumption {\bf A6} holds by an easy inspection of the moduli problem. We let $k$ be an integer and define $\calF=\omg^{\ts k}$ on $X={\bf M}^{\prime}_{H^{\prime}}$. The morphism $\vartheta: w^* \pi^* \calF \ra \calF$ required by {\bf A7} is induced by the pullback morphism via the ($\calO_D$-invariant) projection $\mathcal{A} \ra \mathcal{A}/\mathcal{C}$ where $\mathcal{A},\mathcal{C}$ are the universal abelian scheme and subgroup of order $q$  on $Y$. Assumption {\bf A8} is satisfied because the map $s^*(\pr^* \ts \kappa)$ is induced by pulling back differential forms on $\calA \ts \kappa$ via the Frobenius morphism. For this we are using the fact that $t_1(\Ker({\rm Fr}_q)^{2,1}_1)=\Ker({\rm Fr}_q)$ which follows from condition i) 2) in the definition of the moduli problem and the Cartier duality induced by $\theta$ described above. (See \cite{Kassaei2} \S 10.1.2).

 Now we describe the case of higher levels. We define $X(m)={\bf M}^{\prime}_{H^{\prime},1}(\calP^{m})$, and take $\gerX^0(m)_\rig$ to be the analytification of $N^{\prime}_{H^{\prime},1}(\calP^{m})$. The morphisms $\lambda_{1,\rig}, \lambda_{2,\rig}$ are the analytifications of the morphisms $N^{\prime}_{H^{\prime},1}(\calP^{m}) \ra  M^{\prime}_{H^{\prime},1}(\calP^{m})$ defined, respectively, by forgetting $D$ and dividing by $t_1(D)$. The assumptions {\bf H1} and {\bf H2} are readily satisfied. The morphism $\phi_j:X(m) \ra Y$ is defined as follows. To find $\phi_j(A,i,\theta,Q,\bar{\alpha}^\calP)$, divide  $(A,i,\theta,\bar{\alpha}^\calP)$ by $t_1(<q^{m+1-j}Q>)$ and also enclose the subgroup of order $q$ defined by the image of $q^{m-j}Q$. To be precise,   this definition as it is written works over $\calF_\calP$ but a similar description can be given over $\calO_\calP$ using Drinfeld level structures. Similarly  let $\eta_{j,\rig}$ be the analytification of the morphism $N^{\prime}_{H^{\prime},1}(\calP^{m}) \ra  N^{\prime}_{H^{\prime},0}(\calP)$ defined exactly as $\phi_j$ where one further commands that the extra datum $D$ in the target be generated by the image of its counterpart in the source.  Verifying conditions {\bf H3.1} and {\bf H3.2} are straightforward using Yoneda's lemma. We only make the comment that in verifying {\bf H3.2} one needs observe that if $C$ and $D$ are disjoint $\calO_\calP$-submodule schemes of $A[\varpi]^{2,1}_1$ of rank $q$, then $t_1(C)$ and $t_2(D)$ generate the entire $A[q]$.

We can now apply the results of the paper to this particular example. In particular from Definitions \ref{definition:overconvergence-tame} and \ref{definition:overconvergence-higher} we have spaces of overconvergent modular forms of weight $k \in \ZZ$ over both $M^{\prime}_{H^{\prime},0}(\calP)$  and $M^{\prime}_{H^{\prime},1}(\calP^{m})$. In the first case these spaces are the same as those defined in \cite{Kassaei2}. We can define the $\U=\U_\calP$ operator in both cases using Definition \ref{definition:trace} by taking $c=1/{\rm Nm_{F_\calP/\QQ_p}}(\varpi)$. The definition of the $\U$ operator given in \cite{Kassaei2} differs from this one only in the normalization (there, $c=1/q$ was used). By Propositions \ref{proposition:compcont-tame} and  \ref{proposition:compcont-higher}
we know that $\U_\calP$  is completely continuous on the space of overconvergent modular forms in either case. By Propositions \ref{proposition:analytic-tame} and \ref{proposition:analytic-higher} we find that Buzzard's analytic continuation results hold over these Shimura curves. In other words if $f$ is an overconvergent modular form which is a generalized $\U_\calP$-eigenform with generalized eigenvalue nonzero, then $f$ can be extended to the entire locus where the measure of supersingularity is not equal to $1$. Finally applying Theorem
\ref{theorem:classicality} we have a criterion of classicality for these overconvergent modular forms in all cases (see below). Let us summarize. Recall that $F$ is a totally real field of degree greater than one and $K^\prime$ is small enough as explained at the beginning of this section.

\begin{thm}[Overconvergnce and classicality over unitary Shimura curves in tame and higher levels]
Let $K^\prime$ be either $H^\prime \times K^\prime_0(\calP)$ or $H^\prime \times K^\prime_1(\calP^m)$. We can define spaces of overconvergent modular forms of weight $k \in \ZZ$ on $M^\prime_{K^\prime}$ over which we have the action of a completely continuous operator $\U=\U_\calP$. If an overconvergent modular form $f$ is a generalized eigenform of $\U_\calP$ with  eigenvalue $a_\calP \neq 0$, then it can be extended to the entire nonordinary locus. If   $\val(a_\calP)<k-\val(\rm{Nm}_{F_\calP/\QQ_p}(\varpi))$ then $f$ is classical.

\end{thm}

\subsection{Quaternionic Shimura curves}
Our general results can be   applied to the case of quaternionic Shimura curves as well.   Let $F$ be a totally real field and $\calP$ a prime of $F$ over $p$. Let $B$ be a quaternion algebra over $F$ which is split at $\calP$ and at exactly one infinite place $\tau$. For any choice of an open compact subgroup $K$ of $(B \ts \AA_F^\infty)^\times$ there is a Shimura curve $M_K$ over $F$ whose $\CC$-points are given by
\[
B^\times \backslash (B \ts \AA_F^\infty)^\times \times \gerh^{\pm} /K
\]
where $B^\times$ and $K$ act on $(B \ts \AA_F^\infty)^\times$ by multiplication and on $\CC-\RR=\gerh^\pm$ via, respectively, the inclusion of  $B^\times$ in  $(B \ts_{F,\tau} \RR)^\times \eqra \GL_2(\RR)$ and trivially. When $F=\QQ$ this is a moduli space of abelian surfaces with PEL structure (the so-called {\it fake} or {\it false} elliptic curves). That the set-up of this paper covers this case can be shown in a very similar way to the unitary case using the moduli problem. In this case we take the normalization factor $c$ to be $1/p$. We get the following result:

\begin{thm}[Overconvergnce and classicality over quaternionic Shimura curves / $\QQ$]
Let $f$ be an overconvergent $p$-adic quaternionic modular form of growth condition $r$ with $\val(r)>0$ (where $\val(p)=1$), level $V_1(N)$ with $(p,N)=1$, and weight $k \in \ZZ$ as defined in \S 7 of  \cite{Kassaei1}. If  $f$ is a generalized eigenform of $\U_p$ with eigenvalue $a_p \neq 0$, then it can be extended to the entire nonordinary locus. If furthermore $\val(a_p)<k-1$ then $f$ is classical. Moreover the entire theory can be set up in level $V_1(Np^m)$ for any $m >0$ and the same results hold.
 \end{thm}


\end{document}